\def\no{\if01}
\def\iftwelvept{\no}

\def\ifusepdf{\no}
\def\ifpsfont{\no}
\def\clim{\mathop{\mathrm{colim}}}

\iftwelvept
\documentclass[leqno,12pt]{amsart}
\else
\documentclass[leqno,11pt]{amsart}
\fi
\usepackage{amssymb}
\usepackage{amscd}
\usepackage{latexsym}
\usepackage{verbatim}
\usepackage[all]{xy}

\ifusepdf
\usepackage{hyperref}
\else\fi
\ifpsfont
\usepackage[T1]{fontenc}
\usepackage{times}
\else\fi

\iftwelvept
\else\fi

\theoremstyle{plain}
\newtheorem{Theorem}{Theorem}[section]
\newtheorem{Proposition}[Theorem]{Proposition}
\newtheorem{Lemma}[Theorem]{Lemma}
\newtheorem{Corollary}[Theorem]{Corollary}

\newtheorem{Claim}{Claim}[Theorem]

\theoremstyle{definition}

\newtheorem{Definition}[Theorem]{Definition}
\newtheorem{Remark}[Theorem]{Remark}
\newtheorem{Example}[Theorem]{Example}
\newtheorem*{Conjecture}{Conjecture}

\renewcommand{\theTheorem}{\arabic{section}.\arabic{Theorem}}
\renewcommand{\theClaim}{\arabic{section}.\arabic{Theorem}.\arabic{Claim}}
\renewcommand{\theequation}{\arabic{section}.\arabic{Theorem}.\arabic{Claim}}

\newcommand{\OO}{{\mathcal{O}}}
\newcommand{\LL}{{\mathsf{L}}}

\newcommand{\Hom}{\operatorname{Hom}}
\newcommand{\Ext}{\operatorname{Ext}}

\newcommand{\Image}{\operatorname{Image}}
\newcommand{\Ker}{\operatorname{Ker}}

\newcommand{\Spec}{\operatorname{Spec}}
\newcommand{\Spf}{\operatorname{Spf}}

\newcommand{\SP}{\operatorname{Sp}}

\newcommand{\EXan}{\operatorname{EX}}

\newcommand{\DIM}{\operatorname{dim}}

\newcommand{\an}{\operatorname{rig}}

\newcommand{\lya}{\langle}
\newcommand{\rya}{\rangle}
\newcommand{\lformal}{[\![}
\newcommand{\rformal}{]\!]}

\newcommand{\XX}{\tilde{X}}
\newcommand{\XXX}{\mathcal{X}}
\newcommand{\YYY}{\mathcal{Y}}
\newcommand{\ZZZ}{\mathcal{Z}}

\newcommand{\rig}{\operatorname{rig}}
\newcommand{\aan}{\operatorname{an}}
\newcommand{\Proof}{{\sl Proof.}\quad}
\newcommand{\QED}{{\unskip\nobreak\hfil\penalty50\quad\null\nobreak\hfil
{$\Box$}\parfillskip0pt\finalhyphendemerits0\par\medskip}}

\begin{document}
\title{Deformation theory of rigid-analytic spaces}
\author{Isamu Iwanari}
\begin{abstract}
In this paper, we study a deformation theory of
rigid analytic spaces. We develop a theory of
cotangent complexes for rigid geometry
which fits in with our deformations.
We then use the complexes to give a cohomological
description of infinitesimal deformation of rigid analytic
spaces. Moreover,
we will prove an existence of a formal versal family
for a proper rigid analytic space.
\end{abstract}
\address{Department of Mathematics, Graduate School of Science,
Kyoto University, Kyoto, 606-8502, Japan}
\email{iwanari@math.kyoto-u.ac.jp}
\maketitle

\section*{Introduction}
\renewcommand{\theTheorem}{\Alph{Theorem}}

The purpose of the present paper is to formulate and
develop a deformation theory of rigid analytic
spaces.
We will give a cohomological
description of infinitesimal deformations.
Futhermore we will prove an existence of a formal versal family
for a proper rigid analytic space.

The original idea of deformations
goes back to Kodaira and Spencer.
They developed the theory
of deformations of complex manifolds.
Their deformation theory of complex manifolds is of great
importance and becomes a standard tool for the study of complex analytic spaces.
Our study of deformations of rigid analytic spaces is motivated
by analogy to the case of complex analytic spaces.
More precisely, our interest in the development of such a theory
comes from two sources. First, we want to
construct an analytic moduli theory 
via rigid analytic stacks  by generalizing the classical deformation theory
due to Kodaira-Spencer, Kuranishi and Grauert to the non-Archimedean theory.
This viewpoint will be discussed in Section 5.
Secondly, we may hope that our theory is useful
in arithmetic geometry
no less than the complex-analytic deformation theory is very useful
in the study of complex-analytic spaces.
Actually, our deformation theory will be one of
the key ingredients of the generalization of the theory
of $p$-adic period mappings due to Rapoport-Zink (cf. \cite{RZ}).

Let $K$ be a complete non-Archimedean
valued field and let us consider deformations of affinoid $K$-algebras.
Then we find a serious problem.
As we see later (Remark~\ref{countexample}), it occurs that local deformations of affinoid $K$-algebras
as $K$-algebras are not necessary affinoid algebras.
Thus contrary to the case of schemes, algebraic spaces, etc.,
we can not apply the theory of cotangent complexes of ringed topoi
(cf. \cite{IL}) to our deformation theory.
In this paper, we develop
a general theory of cotangent complexes 
called {\em rigid cotangent complexes} which
fits in with our deformations.
By applying it,  
we will derive
the following theorem.

\begin{Theorem}
\label{intro1}
Let $u:S \rightarrow S'$ be a closed immersion of
rigid analyitic spaces with nilpotent kernel $\mathcal{I}$ of
$u^*:\OO_{S'} \rightarrow \OO_{S}$ such that
$\mathcal{I}^2=0$.
Let $f:X \rightarrow S$ be a flat morphism of rigid
analytic spaces over $K$. 

\begin{enumerate}
\renewcommand{\labelenumi}{(\theenumi)}

  \item The obstruction for the existence of the lifting
 of $f$ to $S'$ lies in
$\Ext^2_{\mathcal{O}_X}(\LL_{X/S}^{\an}, f^* \mathcal{I})$.

  \item If the obstruction $o$ is zero, the set of isomorphism classes of
the liftings of $f$ to $S'$ forms a torsor under
$\Ext^1_{\mathcal{O}_X}(\LL_{X/S}^{\an}, f^*\mathcal{I})$.

  \item Let $\tilde{f}:\tilde{X} \rightarrow S'$ be a flat 
deformation of $X$. Then,
the automorphism group of the lifting $\tilde{X}$ 
is canonically isomorphic to
$\Ext^0_{\mathcal{O}_X}(\LL_{X/S}^{\an}, f^*\mathcal{I})$.
\end{enumerate}
The complex $\LL_{*/*}^{\an}$ in the derived category
$\mathsf{D}^{-}(\mathcal{O}_X \operatorname{-Mod})$ is a rigid cotangent complex
defined in Section $2$.

\end{Theorem}

The idea of rigid cotangent complexes comes from
Raynaud's viewpoint of rigid geometry via formal schemes.
Let $X$ be a quasi-compact and quasi-separated rigid analytic space
over $K$,
and $\XXX$ a formal $R$-schemes which is a formal model of $X$.
Here $R$ is the ring of integer of $K$ (cf. \cite{BL1}).
Then the analytic differential module of $X$ over $K$
is described by the one of $\mathcal{X}$ over $R$ (cf. \cite[1.5]{BL3}).
If $X$ is smooth, the analytic deformations of $X$ is controled by
the analytic diffierencial modules of $X$.
Similarly, the rigid cotangent complex of $X$ over $K$ is constructed
from the cotangent complex of the fixed formal model.
Our theory is closely linked to formal-algebraic geometry
and enable us to deal with deformations
of rigid analytic spaces via formal schemes
and the whole machinery of EGA (for example,
flattening thechniques by blowing-ups are very powerful).
See Section 2 for the details.
We remark that there is another approach to
cotangent complex for rigid geometry due to O. Gabber-L. Ramero.
It uses the theory of Huber's adic spaces 
rather than formal-algebraic geometry.

Our second goal is to show the existence of the formal versal deformation
family which can be viewed as non-Archimedean formal analogue of Kuranishi
and Grauert (cf. \cite{Kura} \cite{Gra}).

\begin{Theorem}
\label{intro2}
Let $X$ be a proper rigid analytic space over $K$.
Then, there exists a formal versal deformation of $X$.
\end{Theorem}

In the appendix, we give a convenient criterion for an
existence of square-zero deformations of a ringed topos.
This criterion may be useful also in the other situations.

{\bf Notations And Conventions}

\smallskip
\begin{enumerate}
\renewcommand{\labelenumi}{(\theenumi)}

 \item Unless otherwise stated, $K$ will be a complete non-Archimedean
valued field.
We denote by $ K \lya X_1,\ldots,X_n \rya$ the
Tate algebra in $n$ indeterminates.

\item For the basic facts and definitions
concerning rigid analytic spaces we refer to \cite{BGR} \cite{BL1}
\cite{BL2} \cite{BL3}.

\item For an affinoid algebra $A$,
we denote by $\SP (A)$ the associated affinoid space (cf. \cite[Chapter 7,8,9]{BGR}).
\item All rigid analytic spaces in this paper will be
quasi-compact and quasi-separated rigid analytic spaces.
Quasi-separatedness means that the diagonal morphism $X \rightarrow X \times X$
is quasi-compact.

\item {\em Properness} means Kiehl's properness (cf. \cite{BGR}).

\item The {\em very weak topology} on an 
affinoid space means Grothendieck topology
such that it has rational subdomains as the admissible open
sets, and its admissible coverings are  unions of  finitely many rational subdomains.
For general rigid analytic spaces, we define it
in the same way. However, unless otherwise stated,
we always equip rigid spaces with the {\em strong topology}
in the sense of \cite{BGR}.

\item By $|\bullet|_{sp}$ we mean the {\em spectral norm} 
(cf. \cite[3.2]{BGR}).

\item Let $\mathcal{A}$ be an admissible formal scheme
(resp. an $\mathcal{O}_{\mathcal{X}}$-module on the admissible formal scheme $\mathcal{X}$,..etc.) in the sense of \cite{BL1}.
Then we denote by $\mathcal{A}^{\rig}$ the rigid analytic spaces
(resp. the $\mathcal{O}_{\mathcal{X}^{\rig}}$-module on the rigid space
$\mathcal{X}^{\rig}$,.. etc.) defined as in \cite[Section 4]{BL1}.
We say that $\mathcal{A}$ is a formal model of $\mathcal{A}^{\rig}$.

\item Let $A$ be a ring and $B$ an $A$-algebra. Let $M$ be a $B$-module.
A short exact
sequence $0 \to M \stackrel{i}{\to} E \stackrel{\pi}{\to} B \to 0$ where $E$
is an $A$-algebra and $\pi$
is a surjective homomorphism of $A$-algebras with $i(M)^2=0$ in $E$ is
said to be a square-zero extension of $B$ by $M$ over $A$.
\end{enumerate}

\subsubsection*{Acknowledgement}
This paper forms a part of my master's thesis. 
I would like to thank Prof. Masao Aoki, Fumiharu Kato and Akio Tamagawa
for helpful conversations and valuable comments concerning the issue presented in this paper, and my adviser Prof. Atsushi Moriwaki
for continuous encouragement. 
I am supported by JSPS
Fellowships for young scientists.

\renewcommand{\theTheorem}{\arabic{section}.\arabic{Theorem}}
\renewcommand{\thesubsubsection}{\arabic{section}.\arabic{subsection}.\arabic{subsubsection}}

\bigskip

\section{First properties of local deformations}

In this section, we will define and prove
first properties of local deformations.

\begin{Definition}
Let
$f:X \rightarrow S$ be a flat morphism
(i.e., a flat morphism of the ringed sites) of rigid analytic spaces
and $S \to S'$ be a closed immersion of
rigid analytic spaces with the nilpotent kernel
$\mathcal{I}:=\Ker (\mathcal{O}_{S'} \to \mathcal{O}_{S})$.
We say that a pair 
$(f': X' \rightarrow S', \phi:X'\times_{S'}S \xrightarrow{\sim} X )$
is a deformation of $f:X \rightarrow S$ to $S'$ if
the following properties are satisfied,
\medskip
\begin{enumerate}
 \item $f'$ is a flat morphism of rigid analytic spaces,

\medskip

 \item $\phi$ is an isomorphism of rigid analytic spaces.

\end{enumerate}

A morphism from 
$(f': X' \rightarrow S', \phi:X'\times_{S'}S \xrightarrow{\sim} X )$ to 
$(f'': X'' \rightarrow S', \psi:X''\times_{S'}S \xrightarrow{\sim} X )$
is a $S'$-morphism $\alpha :X'\to X''$ of rigid analytic spaces
such that $\psi \circ \alpha|_S \circ \phi^{-1}=\operatorname{id}_X$.

\end{Definition}

\begin{Example}
(1) Let $S \to S'$  be a closed immersion of schemes locally
of finite type over $K$ with the nilpotent kernel
$\mathcal{I}:=\Ker(\mathcal{O}_{S'}\to \mathcal{O}_{S})$.
Let $f:X \to S$ be a flat morphism of $K$-schemes and suppose
that $X$ is of finite type over $K$.
Let the pair $(X'/S', X'\times_{S'}S\cong X)$
be a flat deformation of $X$ to $S'$.
Then the analytification $(X^{'\aan}/S^{'\aan}, X^{'\aan}\times_{S^{'\aan}}S^{\aan}\xrightarrow{\sim} X^{\aan})$
is a flat deformation of $X^{\aan}$ to $S^{'\aan}$
.
Here for any scheme $W$ locally of finite type over $K$,
we denote by $W^{\aan}$ the associated rigid analytic space (See \cite{Be} for the analytifications).

(2) Let $T$ be a split $K$-torus and 
$M$ a split lattice of rank $\dim (T)$ in the sense of \cite{BL}.
Note that we can regard $M$ as group $K$-scheme.
If the closed immersion $i:M \to T$
defines $M$ a lattice of full rank in $T$,
the quotient $A:=T/M$ is the rigid analytic group (See \cite[p.661]{BL}).
Let $B$ be an Artin local $K$-algebra. As we see later
(Lemma~\ref{artin:affinoid}), $B$ is an affinoid $K$-algebras.
Let $\tilde{i}:M\times_KB \to T\times_KB$ be a closed immersion
which extends the morphism $i$.
Then the rigid analytic group $(T\times_KB)/(M\times_KB)$
over $\SP (B)$ defines a deformation of $A/K$ to $B$.

\end{Example}

\begin{Remark}
\label{countexample}
It may happen that a square-zero extension of an affinoid algebra is not an affinoid algebra. 
Let $0\rightarrow M\rightarrow E\rightarrow A\rightarrow 0$ be a square-zero extension of an affinoid algebra $A$ over $K$, where $M$ is a finitely generated $A$-module.
By the fundamental theorem due to L. Illusie (cf. \cite[Chapter 3, 1.2.3]{IL}),
the set of isomorphism classes of square-zero extensions of $A$ by
$M$ is classified by the group $\Ext_A^1(\LL_{A/K},M)$ where
$\LL_{A/K}$ is the cotangent complex defined in
\cite[Chapter 2]{IL}.
If $E$ is an affinoid algebra and $A$ is a Tate algebra 
$K \lya T_1, \ldots T_r \rya $, then there exists a splitting $A\rightarrow E$
of $E \to A$ and hence, the extension class of this extension in $\Ext^1_A(\LL_{A/K},M)$ is zero. 
On the other hand, by Theorem~\ref{inj}, there exists a canonical injective
map $\Ext^1_A(\Omega^1_{A/K},M) \to \Ext_A^1(\LL_{A/K},M)$
where $\Omega^1_{A/K}$ is the usual K\"{a}hlar differential module.
Therefore, in order to see the existence of such an extension with $E$ not being an affinoid algebra, it suffices to show the following statement.

\begin{Proposition}
 Let $A=\mathbb{Q}_p\lya T\rya$ be a Tate algebra over $\mathbb{Q}_p$. Then there exists a finitely generated $A$-module $M$ such that $\Ext^1_A(\Omega^1_{A/\mathbb{Q}_p},M)$ is non-zero.
 \end{Proposition}

Proof. First we claim that, for any $A$-module $M$,
we have $\Ext^1_A(\Omega^1_{A/\mathbb{Q}_p},M)=\Ext^1_A(\Ker(\pi),M)$,
where $\Ker (\pi)$ is the kernel of the
homomorphism
\[
0 \to \Ker (\pi ) \to \Omega^1_{A/\mathbb{Q}_p} \stackrel{\pi}{\to} \Omega_{A/\mathbb{Q}_p}^{\rig} \to 0.
\]
Here $\Omega_{A/\mathbb{Q}_p}^{\rig}$ is the differential module defined in
\cite[Section 1]{BL2}.
To prove our claim, we look at the associated long exact sequence to the short exact sequence
\[
\Ext_A^1(\Omega_{A/\mathbb{Q}_p}^{\rig},M) \to \Ext_A^1(\Omega^1_{A/\mathbb{Q}_p},M) \to \Ext_A^1(\Ker (\pi ),M) \to \Ext_A^2(\Omega_{A/\mathbb{Q}_p}^{\rig},M).
\]
Then the claim follows easily from the fact that $\Omega^{\rig}_{A/\mathbb{Q}_p}$ is a free $A$-module.

Thus, what to prove is the existence of a finitely generated $A$-module $M$ such that $\Ext^1_A(\Ker(\pi),M)$ is non-zero.
Suppose we have such an $M$, which is not necessarily finitely generated, then we can actually find a finitely generated $A$-module $M$ having the same property. 
Indeed,
Since $A$ is a principal ideal domain, we have an $A$-injective
resolution, $0\to A \to I \to I' \to 0$.
From the long exact sequence arising from
this resolution, we see that
$\Ext^2_A(\Omega^1_{A/\mathbb{Q}_p},A)=0$.
We have also an $A$-free
resolution $0\to F' \to F \to N \to 0$ because every $A$-submodule of
a free $A$-module is free.
Then we derive $\Ext^1_A(\Omega^1_{A/\mathbb{Q}_p},A)\neq 0$ from the long exact sequence,
\[
\Ext_A^1(\Omega^1_{A/\mathbb{Q}_p},F') \to \Ext_A^1(\Omega^1_{A/\mathbb{Q}_p},F) \to \Ext_A^1(\Omega^1_{A/\mathbb{Q}_p}, N) \to \Ext_A^2(\Omega^1_{A/\mathbb{Q}_p},F')
\]
and $\Ext_A^2(\Omega^1_{A/\mathbb{Q}_p},F')=0$.

The existence of a module $M$ 
with $\Ext^1_A(\Ker(\pi),M) \ne 0$
(not necessarily finitely generated)
follows from the following lemma.

\begin{Lemma}
$\Ker (\pi )$ is a non-zero injective $A$-module.
In particular $\Ker (\pi )$ is not $A$-projective.
\end{Lemma}
\Proof
First, we will show $\Ker (\pi )$ is $A$-injective.
To this aim, since $A$ is
a principal ideal domain, it suffices to prove that
$\Ker (\pi )$ is a divisible $A$-module.
Let $m$ and $a$ be elements in $\Ker (\pi )$
and $A$ respectively.
We want to find an element $n$ of $\Ker (\pi )$ such
that $m=an$.
By induction, we may assume that $a$ is irreducible.
However there exists the exact sequence
\[
0 \to \Hom (\Omega_{A/\mathbb{Q}_p}^{\rig}, A/(a)) \stackrel{\xi}{\to} \Hom (\Omega^1_{A/\mathbb{Q}_p},A/(a))
\to \Hom (\Ker (\pi ),A/(a)) \to 0
\]
where $\xi$ is an isomorphism by 
$\Hom (\Omega_{A/\mathbb{Q}_p}^{\rig}, A/(a)) \cong \operatorname{Der}_{\mathbb{Q}_p}(A,A/(a))
\cong \Hom (\Omega^1_{A/K}, A/(a))$.
Thus $\Ker(\pi )\otimes_{A}A/(a)=0$.
Finally $d(\exp (pT))$  lies in
$\Ker (\pi)$ because
$\exp (pT)=\Sigma_{n=0}^{\infty} (pT)^n/n!$ is
transcendental over $\mathbb{Q}_p(T)$.\QED

\end{Remark}

Next
we investigate the stability of
the Grothendieck topology of a rigid analytic space
under deformations.
First of all, we will consider a class
of morphisms of rigid analytic spaces 
defined as follows.
A morphism $f:X \to Y$ is said to be {\em of affinoid type}
if $f$ is the composite map $X \stackrel{g}{\to} Y\times_K\SP (K \lya T_1,\ldots,T_r \rya ) \stackrel{pr_1}{\to} Y$ for some positive integer $r$ and some closed immersion $g$.

\begin{Proposition}
\label{affinoid:affinoid}
Let $u:S \rightarrow S'$ be a closed immersion of 
rigid analytic spaces with
nilpotent kernel $\Ker (u^*:\mathcal{O}_{S'} \to \mathcal{O}_S)$ and 
X a rigid analytic space. Suppose $f:X \to S$ be a morphism of affinoid type
and $f':\tilde{X} \rightarrow S'$ is a $($possibly non-flat$)$ deformation
of $f$. Then, $f'$ is a morphism of affinoid type.

\end{Proposition}

\Proof
First
we prove the case when
$S:=\SP (A)$, $S':=\SP (A')$ and $I:=\Ker (u^*:A' \to A)$.
Suppose that $I^n=0$.
By induction on $n$, it suffices to
 consider the case where $I^2=0$.
Now we have an exact sequence
\[
\begin{CD}
0 @>>> I\OO_{\tilde{X}} @>>> \OO_{\tilde{X}} @>>> \OO_{X} @>>> 0.
\end{CD}
\]
Since $I^2=0$, $I\OO_{\tilde{X}}=I\OO_{X}$.
Thus by Tate's acyclic theorem (\cite[8.2]{BGR}), we have
\begin{center}
$\check{H}^{1}(\tilde{\mathcal{U}},I\OO_{\tilde{X}})$
$=\check{H}^{1}(\mathcal{U},I\OO_{X})$
$=(0)$
\end{center}
where $\tilde{\mathcal{U}}=\{ \tilde{U}_{i}\}_{i=0}^r$ 
is a finite affinoid cover on $\tilde{X}$ and
$\mathcal{U}=\{ U_i \}_{i=0}^r$
 is the reduction of $\tilde{\mathcal{U}}$ to $X$.
By $\check{H}(\tilde{\mathcal{U}},\bullet)$,
we mean the \v{C}ech cohomology 
with respect to $\tilde{\mathcal{U}}=\{ \tilde{U}_{i}\}$ .
By considering the long exact sequence of the \v{C}ech cohomology with
respect to the very weak topology on $\tilde{X}$, we have an exact sequence
\[
\begin{CD}
0 @>>> H^0 (X,I) @>>> H^0(\tilde{X},\OO_{\tilde{X}})
 @>>> H^0(X,\OO_{X}) @>>> 0.
\end{CD}
\]

We put $H^0(X,\OO_{X})=A \lya X_{1},\ldots,X_{n} \rya /J$.
Let us construct a surjective map from $A'\lya T_1,\ldots T_r\rya$
to $H^0(\tilde{X},\OO_{\tilde{X}})$
which extends a natural surjective map
$A \lya X_{1},\ldots,X_{n} \rya \to H^0(X,\OO_{X})$.
Let 
$\xi_1,\ldots,\xi_n$ be elements of
$H^0(\tilde{X}, \OO_{\tilde{X}})$ which are liftings of $X_{1},\ldots,X_{n},$
respectively. Then, we have the following diagram
\[
\begin{CD}
A' [X_1,\ldots,X_n]   @>>>  A \lya X_1,\ldots X_n \rya   \\
@VV{\pi}V                                     @VVV           \\
H^0(\tilde{X},\OO_{\tilde{X}}) @>>> H^0(X,\OO_X),
\end{CD}
\]
where $\pi$ is defined by $X_{i} \rightarrow \xi _{i}$. 
The right vertical arrow and the lower horizontal arrow are surjective.
We want to see that $\pi$ is continuous
when we equip $K[X_1,\ldots,X_n]$ with the Gaussian norm.
If necessary, we can choose $\xi_1,\dots,\xi_n$
such that 
$\operatorname{max}_{1 \le k \le r}|\operatorname{res}_{\tilde{U}_k}( \xi_i)|_{sp} \le 1$ for $1 \le i \le n$ where $\operatorname{res}_{\tilde{U}_k}( \xi_i)$ is the restriction of $\xi_i$ to $H^0(\tilde{U}_k,\mathcal{O}_{\tilde{X}})$. Indeed, by \cite[3.8.2.2]{BGR},
$|X_i|_{sp} \le |X_i| \le 1$ on $A \lya X_{1},\ldots,X_{n} \rya /J$
and thus we have the same inequality on $U_i$ for $1 \le i \le r$ 
by \cite[3.8.1.4]{BGR}.
Hence, we have $|\operatorname{res}_{\tilde{U}_k}( \xi_i)|_{sp} \le 1$ for any $i,k$.
This implies the composite map 
\begin{center}
$ A'[X_1,\ldots,X_n] \xrightarrow{\pi} $
$ H^0(\tilde{X},\OO_{\tilde{X}})$
$\rightarrow \bigoplus_{i=0}^r H^0(\tilde{U}_i,\OO_{\tilde{U}_i})$
\end{center}
is a continuous map.
Now by Tate's acyclic theorem we have
\begin{center}
$H^0(\tilde{X},\OO_{\tilde{X}})=\Ker (\bigoplus_{i=0}^r H^0(\tilde{U}_i,
\OO_{\tilde{U}_i})$
$\rightrightarrows \bigoplus_{i<j}H^0(\tilde{U}_{i}\cap \tilde{U}_{j},
 \OO_{\tilde{U}_{i}\cap \tilde{U}_{j}}))$.

\end{center}
Note that the topology of $\bigoplus_{i=0}^r H^0(\tilde{U}_i,
\OO_{\tilde{U}_i})$
(i.e. direct sum of the topologies on affinoid algebras
 $H^0(\tilde{U}_i,\OO_{\tilde{U}_i}))$
induces the topology on $H^0(\tilde{X},\OO_{\tilde{X}})$
which is complete. 
Thus we see $\pi$ is also continuous.
From this, there exists a unique homomorphism
\begin{center}
$\Pi: A' \lya X_1,\ldots,X_n \rya \longrightarrow $
$ H^0(\tilde{X},\OO_{\tilde{X}})$
\end{center}
which extends the homomorphism $\pi$.
We claim that $\Pi$ is surjective. Indeed, the kernel of the surjection
$H^0(\tilde{X},\OO_{\tilde{X}}) \rightarrow H^0(X,\OO_X)$ is
$IH^0(\tilde{X},\OO_{\tilde{X}})$. 
On the other hand, it is clear that the image of $\pi$
contains $IH^0(\tilde{X},\OO_{\tilde{X}})$.
Hence, we see that $\Pi$ is a surjection.
Therefore, there exists the commutative diagram of rigid analytic spaces
\[
\begin{CD}
X                       @>>{i}> \tilde{X}           \\
@V{\alpha}VV                           @VV{\beta}V  \\
\mathbb{D}_{A}^{n}      @>>{j}> \mathbb{D}_{A'}^{n}
\end{CD}
\]
where $\mathbb{D}_{A}^{n}$ and $\mathbb{D}_{A'}^{n}$
are $\SP (A \lya X_1,\ldots,X_n \rya )$ 
and $\SP (A' \lya X_1,\ldots,X_n \rya)$
respectively. $i$ and $j$ are maps which induce
bijective maps of underlying sets.
The morphism $j \circ \alpha$
 induce an injective map of undelying sets.
Hence $\beta$ is bijective onto the
image of $j \circ  \alpha$ as sets. To prove the proposition, it suffices to
show that homomorphism $\OO_{\mathbb{D}_{A'}^n,\beta (x)} \xrightarrow{\beta^*}$
$\OO_{\tilde{X},x}$ is surjection for any point $x$ of $\tilde{X}$.
Let $a$ be an element of ${\OO}_{\tilde{X},x}$.
There is an element $b$ in ${\OO}_{\mathbb{D}_{A}^n,\beta (x)}$
such that ${\alpha}^*(b)={i}^*(a)$.
Moreover there is an element 
$c$ in ${\OO}_{\mathbb{D}_{A'}^n,\beta (x)}$
such that ${j}^*(c)=b$.
Then, ${\beta}^*(c)-a$ is in 
$\Ker ({i}^*)=I {\OO}_{\tilde{X},x}$.
Since $I^2=0$, we have 
$I\Image ({\beta}^*)=I{\OO}_{\tilde{X},x}$.
Thus we have an element $d$ in ${\OO}_{\mathbb{D}_{A'}^n,\beta (x)}$
such that ${\beta}^*(d+c)=a$, since $I$ is contained by
$\Image ({\beta}^*)$. 

Finally the assertion for the general case follows from the proof of the local case.
\QED

\begin{Proposition}
\label{rat:rat}
Let $u: A'\rightarrow A$ be a surjective homomorphism of 
affinoid $K$-algebras
with the nilpotent kernel and $X$ be a $K$-affiniod space
over $\SP (A)$. Let $U$ be a rational subdomain of $X$.
 Suppose $\tilde{X} \rightarrow \SP (A')$ is a $($possibly non-flat$)$ deformation
of $X$ and $\tilde{X}$ is an affinoid spaces. 
Then, the lifting $\tilde{U}$ of $U$ in $\tilde{X}$ is a rational
subdomain of $\tilde{X}$.

\end{Proposition}

\Proof
We put 
\begin{align*}
U &=X(f_0,\dots,f_r)  \\
  &=\{ x \in X | |f_1(x)|\le |f_0(x)|,\dots,|f_r(x)| \le
|f_0(x)| \} \\
  &=\SP (H^0(X, \OO_{X}) \lya T_1,\dots,T_r \rya/(f_1-T_1f_0,\dots,f_r-T_rf_0))  
\end{align*}
where $f_0,\dots,f_r$ generates the unit ideal.
We choose elements $\tilde{f}_0,\dots,\tilde{f}_r \in H^0(\tilde{X},
\OO_{\tilde{X}})$ which are the liftings of $f_0,\dots,f_r$
respectively. Since $X \rightarrow \tilde{X}$ is a nilpotent thickening,
we have
\begin{align*}
\tilde{U}&=\tilde{X}(\tilde{f}_0,\dots,\tilde{f}_r)  \\
         &=\{ x \in \tilde{X} | |\tilde{f}_1(x)|\le |\tilde{f}_0(x)|,
\dots,|\tilde{f}_r(x)| \le|\tilde{f}_0(x)| \} \\
         &=\SP (H^0(\tilde{X}, \OO_{\tilde{X}}) \lya T_1,\dots,T_r \rya /
(\tilde{f}_1-T_1\tilde{f}_0,\dots,\tilde{f}_r-T_r\tilde{f}_0)).
\end{align*}
This implies the proposition.
\QED

\begin{Theorem}
\label{adm:adm}
Let $u:S \rightarrow S'$ be a closed immersion of 
rigid analytic spaces with
nilpotent kernel
$\Ker (u^*:\mathcal{O}_{S'} \to \mathcal{O}_S)$ and $f:X \to S$
a morphism of rigid analytic spaces.
Let $f':\tilde{X}\to S'$ be a $($possibly non-flat$)$ deformation of $f$.
Suppose $U$ is an admissible open set with respect
to the strong topology on $X$.
Then, the lifting $\tilde{U}$ of $U$ in $\tilde{X}$ is
the admissible open set of $\tilde{X}$.
Similarly, an admissible covering of $X$ lifts to admissible covering of
$\tilde{X}$. In particular the Grothendieck topology of a rigid analytic
space is stable under nilpotent deformations.
\end{Theorem}

\Proof
By Proposition~\ref{affinoid:affinoid} and \cite[9.1.2.2]{BGR} ,
 we may assume that $X$, $\tilde{X}$, $S$ and $S'$ are  affinoid spaces.
Next, we note that the strongest topology among
the topologies which are slightly finer (See for the definition
\cite[9.1.2.1]{BGR}) than very weak topology
coincides with the strong topology by the theorem of Gerritzen-Grauert
\cite[7.3.5.3]{BGR}.

Thus to prove the assertion,
it suffices to check that (cf. \cite[9.1.4.2]{BGR}):

\medskip

\begin{enumerate}

 \item The set $\tilde{U}$ admits a covering $\{ \tilde{U}_i\}_i$
by affinoid subdomains $\tilde{U}_i\subset \tilde{X}$ such that, for any
affinoid morphism 
$\phi :Y \rightarrow \tilde{X}$ with $\phi (Y) \subset \tilde{U}$,
the covering $\{ \phi^{-1}(\tilde{U}_i) \}_i$ of $Y$ has a finite 
rational subdomain covering
which refines it.

\medskip

 \item  Let $\{ \tilde{V}_j \}_j$ be a covering
which are the liftings
of an admissible covering of the admissible open set $V \subset X$.
Note that if the first half of our claim is verified, $V_i$ are admissible.
Let $\tilde{V} \subset \tilde{X}$ be the lifting of $V$.
Then, for any affinoid morphism $\phi :Y \rightarrow \tilde{X}$
with $\phi (Y) \subset \tilde{V}$, the covering $\{ \phi^{-1}(\tilde{V}_i) \}_i$
of $Y$ has a finite rational subdomain covering which refines it.

\end{enumerate}

Since $|X|=|\tilde{X}|$, we can easily check these conditions
by using Proposition~\ref{rat:rat}.
\QED

\section{Rigid cotangent complex}

Since the theory of cotangent complex for ringed topoi (cf. \cite{IL})
can not be applied to our problem, we define a {\em rigid
cotangent complex} inspired by \cite[Proposition 1.5]{BL3}.

\subsection{Analytic cotangent complex of formal schemes}
First of all, following \cite{GR}, we will recall
the analytic cotangent complexes for formal schemes
locally of topologically finite presentation over $\Spf R$,
where $R$ is the ring of integer of 
$K$. Let $\pi$ be an element in the maximal ideal of $R$
and $A$ a complete $R$-algebra topologically of finite presentation.
Consider the $\pi$-adic completion functor
\[
(A\operatorname{-Mod})\longrightarrow (A\operatorname{-Mod})
\]
\[
L \mapsto L^{\wedge}.
\]
This functor induces the derived functor
\[
\mathsf{D}^{-}(A\operatorname{-Mod}) \longrightarrow \mathsf{D}^{-}(A\operatorname{-Mod}).
\]
Indeed, 
for any quasi-isomorphism of complexes of flat $A$-modules
$K^{\bullet}\to L^{\bullet}$,
the induced homomorphism
$(K^{\bullet})^{\wedge} \to (L^{\bullet})^{\wedge}$
($(\operatorname{complex})^{\wedge}$ denotes the termwise completion)
is a quasi-isomorphism (cf. \cite[7.1.11]{GR}).
Moreover, the derived category $\mathsf{D}^{-}(R\operatorname{-Mod})$
is naturally identified with 
the localization of 
the homotopy category $K^-(A\operatorname{-flat\ Mod})$
up to quasi-isomorphisms. Thus
$\pi$-adic completion functor induces the derived functor.

Let $\phi :A\rightarrow B$ be a homomorphism of complete $R$-algebras
of topologically finite presentation.
The $B$-module of analytic differentials relative to $\phi$
is defined as $\Omega_{B/A}^{\aan}:=\Omega_{B/A}^{1\wedge}$.
The analytic cotangent complex of $\phi$ is the complex
$\LL_{B/A}^{\aan}:=\LL_{B/A}^{\wedge}$.
Here $\LL_{B/A}$ is the usual cotangent complex of $\phi$ (cf. \cite{IL}).
There exists a natural map $\LL_{B/A}^{\aan} \rightarrow \Omega_{B/A}^{\aan}$
and an isomorphism $H_0(\LL_{B/A}^{\aan})\cong \Omega_{B/A}^{\aan}$.
If $\phi$ is smooth, then there is a natural quasi-isomorphism
$\LL_{B/A}^{\aan} \cong \Omega_{B/A}^{\aan}[0]$.

Next we define the analytic cotangent complex 
for formal schemes locally of topologically finite presentation
over $\Spf R$ by gluing the complexes
constructed as above.
Let $f:\mathcal{X}\rightarrow \mathcal{Y}$ be a morphism of
formal schemes locally of topologically finite presentation over $\Spf R$ (cf. \cite{BL1}).
We suppose that $\mathcal{Y}$ is affine for a while. For an affine open set $\mathcal{U}$ in $\mathcal{X}$, the small category $F_{\mathcal{U}}$
of all affine open sets $\mathcal{V}$ in $\mathcal{Y}$ 
with $f(\mathcal{U}) \subset \mathcal{V}$ is
a cofiltered family under inclusion.
For every $\mathcal{V} \in F_{\mathcal{U}}$, $\mathcal{O}_{\mathcal{Y}}(\mathcal{V})$ is
a complete $R$-algebra of topologically finite presentation. The cofiltered family of maps
$\mathcal{O}_{\mathcal{Y}}(\mathcal{V})\rightarrow \mathcal{O}_{\mathcal{X}}(\mathcal{U})$ gives rise to the correspondence
\[
\mathcal{U} \longrightarrow \LL (\mathcal{U}/\mathcal{Y}):=\displaystyle
\clim_{\mathcal{V} \in F_{\mathcal{U}}} \LL^{\aan}_{\mathcal{O}_{\mathcal{X}}(\mathcal{U})/\mathcal{O}_{\mathcal{Y}}(\mathcal{V})}.
\]
Now note that every usual cotangent complex is constructed
as a complex of free modules in a functorial fashion.
Hence for any homomorphism $A\to B$, we can construct
the analytic cotangent complex of $A \to B$ as a complex of flat $B$-modules
in a functorial fashion.
(The flatness follows from the fact that
for any $B$-flat module $F$, the completion $F^{\wedge}$ is
also $B$-flat (cf. for example \cite[7.1.6 (1)]{GR})).
Therefore,
by \cite[Chapter 0, 3.2.1]{EGA1}, we can extend the above correspondence
to the complex of presheaves on $\mathcal{X}$.
The analytic cotangent complex is defined as the complex of termwise the
associated sheaves. Finally for a general
 formal scheme $\mathcal{Y}$, we can construct the
complex by gluing
the complex constructed locally on $\mathcal{Y}$.

The following proposition will be used in the next subsection.

\begin{Proposition}
\label{formal-base-change}
Let $f:\XXX \to \YYY$ and $g:\ZZZ \to \YYY$ be morphisms of
formal schemes locally of topologically finite presentation over $\Spf R$.
Let $\XXX \times_{\YYY}\ZZZ$ be the fibre product
of $\XXX$ and $\ZZZ$ over $\YYY$ in the category
of formal $R$-schemes.
 Then there exists a natural quasi-isomorphism
\[
L\operatorname{pr}_1^*\LL_{\XXX /\YYY}^{\aan} \stackrel{\sim}{\to} \LL_{\XXX\times_{\YYY}\ZZZ /\ZZZ}^{\aan}
\]
where $\operatorname{pr}_1:\XXX \times_{\YYY}\ZZZ \to \XXX$ is the first projection.

\end{Proposition}
\Proof
Since our assertion is a local issue,
we may suppose that $\XXX$, $\YYY$ and $\ZZZ$ are affine.
Set $\XXX=\Spf B$, $\YYY=\Spf A$ and $\ZZZ=\Spf C$.
First we shall show that there exists a natural quasi-isomorphism
$\phi:(B\hat{\otimes}_AC)\otimes_{B}\LL_{B/A}^{\wedge} \to
\LL_{(B\hat{\otimes}_AC)/C}^{\wedge}$.
We remark that $\Spf B\hat{\otimes}_AC$ is the fibre product of $\Spf B$ and
$\Spf C$
over $\Spf A$ in the category of formal $R$-schemes.
In this proof, unless otherwise stated, we view complexes
as just complexs (not objects in the derived category)
and we denote by $\stackrel{\operatorname{qis}}{\to}$ a quasi-isomorphism
of complexes.
Actually, usual cotangent complexed are constructed as complexes
of flat modules via
standard resolutions and thier completions are also flat.
There exists natural quasi-morphisms
$((B\hat{\otimes}_AC)\otimes_{(B\otimes_AC)}(B\otimes_AC)\otimes_B\LL_{B/A})^{\wedge}\stackrel{\sim}{\to}((B\hat{\otimes}_AC)\otimes_B\LL_{B/A})^{\wedge}\stackrel{\operatorname{qis}}{\to}((B\hat{\otimes}_AC)\otimes_B\LL_{B/A}^{\wedge})^{\wedge}$.
The second quasi-isomorphism follows from \cite[Lemma 7.1.22]{GR}.
On the other hand, by the base change theorem of usual cotangent complexes
\cite[Chapter 2, 2.2]{IL}, we have a natural quasi-isomorphism
$((B\hat{\otimes}_AC)\otimes_{(B\otimes_AC)}(B\otimes_AC)\otimes_B\LL_{B/A})^{\wedge} \stackrel{\operatorname{qis}}{\to}((B\hat{\otimes}_AC)\otimes_{(B\otimes_AC)}\LL_{(B\otimes_AC)/C})^{\wedge}$.
The next claim implies a quasi-isomorphism $(B\hat{\otimes}_AC)\otimes_B\LL_{B/A}^{\wedge}\stackrel{\operatorname{qis}}{\to}((B\hat{\otimes}_AC)\otimes_B\LL_{B/A}^{\wedge})^{\wedge}$
and thus we have a natural quasi-isomorphism
$(B\hat{\otimes}_AC)\otimes_B\LL_{B/A}^{\wedge}\stackrel{\sim}{\to}
((B\hat{\otimes}_AC)\otimes_{(B\otimes_AC)}\LL_{(B\otimes_AC)/C})^{\wedge}$
in the derived category.

\begin{Claim}
\label{help}
Let $A\to B$ be a homomorphism of admissible $R$-algebras and
$K^{\bullet}:=\LL^{\aan}_{B/A}$ the analytic cotangent complex.
Then there is a natural quasi-isomorphism
\[
K^{\bullet} \stackrel{\LL}{\otimes}_AB \stackrel{\sim}{\to} (K^{\bullet}\stackrel{\LL}{\otimes}_AB)^{\wedge}
\]
in the derived category $\mathsf{D}^{-}(A\operatorname{-Mod})$.
\end{Claim}
\Proof
It suffices to show that for any positive integer $n$,
the truncation
\[
\tau_{[-n}(K^{\bullet} \stackrel{\LL}{\otimes}_AB) \stackrel{\sim}{\to} \tau_{[-n}(K^{\bullet}\stackrel{\LL}{\otimes}_AB)^{\wedge}
\]
is a quasi-isomorphism.
Due to \cite[7.1.15, 7.1.33 (1)]{GR} and the fact that $B$ is coherent, we can assume that $\tau_{[-n}K^{\bullet}$ is a complex of free $A$-modules of finite type.
Therefore it suffices to show that for every free $A$-module of finite type
$F$, the natural homomorphism $F\otimes_AB \to (F\otimes_AB)^{\wedge}$ is an isomorphism.
However it is clear.
\QED

Therefore, to see that $\phi$ is a quasi-isomorphism,
it suffices to show that the natural morphism
$\psi:((B\hat{\otimes}_AC)\otimes_{B\otimes_AC}\LL_{(B\otimes_AC)/C})^{\wedge}
\to \LL_{(B\hat{\otimes}_AC)/C}^{\wedge}$ is
a quasi-isomorphism.
Note that, by the base-chage theorem
of usual cotangent complexes \cite[Chapter 2, 2.2]{IL}, 
$((B\hat{\otimes}_AC)\otimes_{B\otimes_AC}\LL_{(B\otimes_AC)/C})^{\wedge}_n
:=((B\hat{\otimes}_AC)\otimes_{B\otimes_AC}\LL_{(B\otimes_AC)/C})^{\wedge}\otimes_R(R/\pi^nR)\stackrel{\operatorname{qis}}{\to} \LL_{(B_n\otimes_{A_n}C_n)/C_n}$,
$(\LL_{(B\hat{\otimes}_AC)/C}^{\wedge})_n:=\LL_{(B\hat{\otimes}_AC)/C}^{\wedge}
\otimes_R(R/\pi^nR) \stackrel{\operatorname{qis}}{\to} \LL_{(B_n\otimes_{A_n}C_n)/C_n}$
where $B_n=B\otimes_R(R/\pi^nR)$, $A_n=A\otimes_R(R/\pi^nR)$ and
$C_n=C\otimes_R(R/\pi^nR)$.
Thus the reduction $\psi\otimes_R(R/\pi^nR)$ is a quasi-isomorphism.
Now consider the right derived functor
\[
R \lim :\mathsf{D}((R\operatorname{-Mod})^{\mathbb{N}})) \to
\mathsf{D}(R\operatorname{-Mod})
\]
where $(R\operatorname{-Mod})^{\mathbb{N}}$ is the projective
system of $R$-modules and $\lim :(R\operatorname{-Mod})^{\oplus \mathbb{N}}$
$\to R\operatorname{-Mod}, \{ M_i, d_i:M_i \to M_{i-1}\}_{i\ge 1} \mapsto \operatorname{proj.lim}_iM_i$
is the inverse limit functor.
The projective systems of complexes
$((B\hat{\otimes}_AC)\otimes_{B\otimes_AC}\LL_{(B\otimes_AC)/C})^{\wedge}_n)_{n \ge 1}$
and $((\LL_{(B\hat{\otimes}_AC)/C}^{\wedge})_n)_{n\ge 1}$ are acyclic
for the fucntor $\lim$ because they are concisting of surjections.
Thus we have that
\[
R \lim ((B\hat{\otimes}_AC)\otimes_{B\otimes_AC}\LL_{(B\otimes_AC)/C})^{\wedge}_n)_{n \ge 1}=((B\hat{\otimes}_AC)\otimes_{B\otimes_AC}\LL_{(B\otimes_AC)/C})^{\wedge}
\]
and
\[
R \lim ((\LL_{(B\hat{\otimes}_AC)/C}^{\wedge})_n)_{n\ge 1}=\LL_{(B\hat{\otimes}_AC)/C}^{\wedge}.
\]
Since there exists a quasi-isomorphism
$R \lim ((B\hat{\otimes}_AC)\otimes_{B\otimes_AC}\LL_{(B\otimes_AC)/C})^{\wedge}_n)_{n \ge 1}\stackrel{\operatorname{qis}}{\to}  R \lim ((\LL_{(B\hat{\otimes}_AC)/C}^{\wedge})_n)_{n\ge 1}$, we see that $\psi$ is a quasi-isomorphism and thus $\phi$ is a quasi-isomorphism.
Futhermore, a natural homomorphism $(B\hat{\otimes}_AC)\stackrel{\LL}{\otimes}_{B}\LL_{B/A}^{\wedge} \to
\LL_{(B\hat{\otimes}_AC)/C}^{\wedge}$ is a quasi-isomorphism (in the derived category) because $\LL_{B/A}^{\wedge}$ consists of flat $B$-modules.

Next note that
$(B \hat{\otimes}_{A}C)\stackrel{\LL}{\otimes}_{B}\LL_{B/A}^{\wedge}$ and
$\LL_{(B\hat{\otimes}_AC)/C}^{\wedge}$ are pseudo-coherent
\footnote{Let $n$ be a integer. We say that a complex of $R$-module
$K^{\bullet}$ is $n$-pseudo-coherent if there exists
a quasi isomorphism $C^{\bullet}\to K^{\bullet}$
where $C^{\bullet}$ is a complex bounded above and
$C^k$ is a finitely generated free $R$-module for every $k \ge n$.
We say that $K^{\bullet}$ is pseudo-coherent if
$K^{\bullet}$ is $n$-pseudo-coherent for every integer $n$.
$n$-pseudo-coherence is stable under quasi-isomorphisms.}.
Let $L_i^{\Delta}$ (resp. $M_i^{\Delta}$) be the sheaf of the coherent
$\mathcal{O}_{\Spf (B \hat{\otimes}_{A}C)}$-module
associated to the
coherent $B \hat{\otimes}_{A}C$-module
$L_i:=H_i((B \hat{\otimes}_{A}C)\stackrel{\LL}{\otimes}_{B}\LL_{B/A}^{\wedge})$
(resp. $M_i:=H_i(\LL_{(B\hat{\otimes}_AC)/C}^{\wedge})$).
To complete the proof of the proposition,
it suffices to show that following claim.

\begin{Claim}
Under the same assumption as above,
there exists natural isomorphisms
\[
L_i^{\Delta}\cong H_i(L\operatorname{pr_1}^*\LL_{\Spf B/\Spf A}^{\aan})
\]
and
\[
M_i^{\Delta}\cong H_i(\LL_{(\Spf B\hat{\otimes}_AC)/\Spf C}^{\aan}).
\]
\end{Claim}
\Proof
The second assertion can be shown by the same way of the proof
of the first assertion. Hence we will prove the first
quasi-isomorphism.
What to prove is that the natural isomorphism
$\xi:L_i^{\Delta}\to H_i(L\operatorname{pr_1}^*\LL_{\Spf B/\Spf A}^{\aan})$
induces an isomorphism on each stalk.
Let $U:=\Spf R(U)$, $V:=\Spf R(V)$, and $W:=\Spf R(W)$ be affine open sets of $\Spf  B\hat{\otimes}_AC$,
$\Spf B$, and $\Spf A$ respectively and suppose that
$\operatorname{pr_1}(U)\subset V$ and $f(V)\subset W$.
By the transitivity of analytic cotangent complexes
\cite[7.1.33 (2)]{GR}, we easily see that
there exists a natural isomorphism
$H_i(R(U)\otimes_{R(V)}\LL^{\wedge}_{R(V)/R(W)})\cong
H_i((R(U)\otimes_{B}\LL^{\wedge}_{B/A})$.
Therefore we have natural isomorphisms
\begin{align*}
H_i(L\operatorname{pr_1}^*\LL_{\Spf B/\Spf A}^{\aan})_x
&\cong \displaystyle\clim_{x \in U \atop \operatorname{affine\ open}}\displaystyle\clim_{\operatorname{pr_1}(U) \subset V \atop
\operatorname{affine\ open}}\displaystyle\clim_{f(V)\subset W \atop
\operatorname{affine\ open}}H_i(R(U)\otimes_{R(V)}\LL_{R(V)/R(W)}^{\wedge}) \\
&\cong \displaystyle\clim_{x \in U \atop \operatorname{affine\ open}}
H_i(R(U)\otimes_{B}\LL_{B/A}^{\wedge}) \\
&\cong \displaystyle\clim_{x \in D(c) \atop c \in B\hat{\otimes}_AC}
H_i(R(D(c))\otimes_{B}\LL_{B/A}^{\wedge}) \\
&\cong L_i\otimes_{(B\hat{\otimes}_AC)}\OO_{\Spf (B\hat{\otimes}_AC),x}
\end{align*}
where $D(c)=\{ x \in \Spf (B\hat{\otimes}_AC),c \notin \mathsf{m}_x \} $
and $R(D(c)):=H^0(D(c),\mathcal{O}_{\Spf (B\hat{\otimes}_AC)})$.
The final isomorphism follows from the next lemma.

\begin{Lemma}
Let $A$ be a complete $R$-algebra of topologically finite presentation.
For any $c \in A$, let $D(c):=\{ x\in \Spf A;c \notin \mathfrak{m}_x\}$.
Then the natural map $A \to H^0(D(c),\OO_{\Spf A})$ is flat.

\end{Lemma}
\Proof
Note that $H^0(D(c),\OO_{\Spf A})$ is the $\pi$-adic completion of $A_c$.
Then our claim follows from the fact that the $\pi$-adic completion of every flat A-module
is also $A$-flat.
\QED

\subsection{Rigid cotangent complex}
 In the rest of this section we often use
the terminologies in \cite{BL1} \cite{BL2} \cite{BL3}.
For the definition of admissible formal schemes,
admissible blow-up, smoothness, etc., we refer to
\cite{BL1} \cite{BL2} \cite{BL3}.
Before we define rigid cotangent complexes,
we will show some results which are used later.

\begin{Lemma}
\label{defined}
\begin{enumerate}
\renewcommand{\labelenumi}{(\theenumi)}

 \item Let $ X \stackrel{f}{\to} Y \stackrel{g}{\to} Z$ be a sequence of morphisms of rigid analytic
spaces. Then there is a natural quasi-isomorphism
\[
L(g \circ f)^* \stackrel{\sim}{\to} Lf^* \circ Lg^*.
\]

\item Let $\mathcal{X} \stackrel{\tilde{f}}{\to} \mathcal{Y}
\stackrel{\tilde{g}}{\to} \mathcal{Z}$ be a
sequence of admissible formal $R$-schemes which is a
formal model of $ X \stackrel{f}{\to} Y \stackrel{g}{\to} Z$
$($Actually we can choose such a sequence of formal schemes
for any sequence of rigid analytic spaces by \cite[Theorem 4.1]{BL1}$)$.
Then there is a natural quasi-isomorphism
\[
(L \tilde{f}^*(\LL_{\mathcal{Y}/\mathcal{Z}}^{\aan}))^{\rig}
\stackrel{\sim}{\to}
Lf^{*}(\LL_{\mathcal{Y}/\mathcal{Z}}^{\aan})^{\rig}.
\]

\end{enumerate}
\end{Lemma}

\Proof
(1) Our claim follows from Grothendieck spectral sequence (cf. \cite[1.8.7]{KS}).

(2) First, note that the assertion is local on $\mathcal{X}$.
Thus we suppose that $\mathcal{X}$ and $\mathcal{Y}$ are affine.
Put
$A=H^0(\mathcal{X},\mathcal{O}_{\mathcal{X}})$ and
$B=H^0(\mathcal{Y},\mathcal{O}_{\mathcal{Y}})$.
Moreover, we may replace $\LL_{\mathcal{Y}/\mathcal{Z}}^{\aan}$
by $\LL(n):=\tau_{[-n}\LL_{\mathcal{Y}/\mathcal{Z}}^{\aan}$
and prove our assertion for the latter complexes for every integer $n$.
Take a complex of $B$-flat modules of finite presented
$\LL_0^{\bullet}$ which represents the complex $\LL(n)$.
Then we have
$(L \tilde{f}^*(\LL(n))^{\rig}
=(\LL_0^{\bullet} \otimes_BA \otimes_{R}K)^{\sim}$
(By $(\bullet)^{\sim}$ we denote the associated sheaf). On the other hand,
$\LL_0^{\bullet}\otimes_RK$ is a complex of $B\otimes_{R}K$-flat module 
which represents $\LL(n)^{\rig}$. Thus we have a quasi-isomorphism
$Lf^*\LL(n) \cong
((\LL_0^{\bullet}\otimes_RK)\otimes_{(B\otimes K)}(A\otimes K))^{\sim}$.
Therefore, we have a natural quasi-isomorphism
$(\LL(n)\stackrel{\LL}{\otimes}_BA)^{\rig}
\stackrel{\sim}{\to}
\LL(n)^{\rig}\stackrel{\LL}{\otimes}_{(B\otimes_RK)}(A\otimes_RK)$.
\QED

Let $f:X \rightarrow Y$ be a morphism of rigid analytic spaces over $K$.
There exists
a morphism of admissible formal $R$-schemes
\[
\tilde{f}:\mathcal{X} \longrightarrow \mathcal{Y}
\]
which is a formal model of $f$ (cf. \cite[Theorem 4.1]{BL1}).
\begin{Proposition}
\label{core}
The complex $(\LL^{\aan}_{\mathcal{X}/\mathcal{Y}})^{\an}$ is independent of
the choice of the formal model
$\tilde{f}\colon \mathcal{X} \rightarrow \mathcal{Y}$ of $f\colon X \rightarrow Y$, i.e., depends
only on the morphism $f:X \rightarrow Y$.
\end{Proposition}

\Proof Let $\tilde{f}':\mathcal{X}' \rightarrow \mathcal{Y}'$ be another formal
model of $f:X \rightarrow Y$. By an easy application of the theorem of
Raynaud \cite[Theorem 4.1]{BL1}, we can find the following commutative
diagram
\[
 \xymatrix{
 \mathcal{X} \ar[d]^{\tilde{f}} & \mathcal{X}'' \ar[l]_{\alpha} \ar[r]^{s} \ar[d]^{\tilde{f}''} & \mathcal{X}' \ar[d]^{\tilde{f}'}  \\
 \mathcal{Y} & \mathcal{Y}'' \ar[l]_{\beta} \ar[r]^{t} &  \mathcal{Y}  \\
  }
\]
such that $\tilde{f}, \tilde{f}', \tilde{f}''$ are formal models
of $f:X \rightarrow Y$ and $\alpha, \beta, s, t$ are admissible blowing
ups. It suffices to show  that there are natural isomorphisms
$(\LL_{\mathcal{X}/\mathcal{Y}}^{\aan})^{\rig} \cong
(\LL_{\mathcal{X}''/\mathcal{Y}''}^{\aan})^{\rig} \cong
(\LL_{\mathcal{X}'/\mathcal{Y}'}^{\aan})^{\rig}$.
We will show an existence of the first isomorphism.
The second one follows from
the proof of the first one.
To this end, consider the commutative diagram
\[
 \xymatrix{
 \mathcal{X}'' \ar[rdd]_{\alpha} \ar[rd]_{\eta} \ar[rrd]^{\tilde{f}''} &  &  \\
  &\mathcal{X}\times_{\mathcal{Y}}\mathcal{Y}'' \ar[r]_{pr_2} \ar[d]^{pr_1} & \mathcal{Y}''\ar[d]^{\beta} \\
    & \mathcal{X} \ar[r]^{f} & \mathcal{Y}  \\
  }
\]
where $\mathcal{X}\times_{\mathcal{Y}}\mathcal{Y}''$ is the fiber product of
$\mathcal{X}$ and $\mathcal{Y}''$ over $\mathcal{Y}$ in the category of formal $R$-schemes and $\eta$ is the induced morphism by $\alpha$ and $\tilde{f}''$.

By the transitivity to the sequence $\mathcal{X}'' \rightarrow \mathcal{X}\times_{\mathcal{Y}}\mathcal{Y}'' \rightarrow \mathcal{Y}''$ (\cite[7.2.13]{GR}), there exists
a distinguished triangle
\[
L\eta^{*} \LL_{\mathcal{X}\times_{\mathcal{Y}}\mathcal{Y}''/\mathcal{Y}''}^{\aan}
\rightarrow \LL_{\mathcal{X}''/\mathcal{Y}''}^{\aan} \rightarrow
\LL_{\mathcal{X}''/\mathcal{X}\times_{\mathcal{Y}}\mathcal{Y}''}^{\aan} \rightarrow
L\eta^{*} \LL_{\mathcal{X}\times_{\mathcal{Y}}\mathcal{Y}''/\mathcal{Y}''}^{\aan}[1].
\]
Now consider the functor of derived categories
\[
\operatorname{Rig}:\mathsf{D}^{-}(\mathcal{O}_{\mathcal{X}''}\operatorname{-Mod})\rightarrow
\mathsf{D}^{-}(\mathcal{O}_{X}\operatorname{-Mod}),
\]
which is induced by $M \mapsto (M)^{\rig}$ for any $\mathcal{O}''$-module $M$.
By applying this functor, we obtain
the following distinguished triangle,
\[
(L\eta^{*} \LL_{\mathcal{X}\times_{\mathcal{Y}}\mathcal{Y}''/\mathcal{Y}''}^{\aan})^{\rig}
\rightarrow (\LL_{\mathcal{X}''/\mathcal{Y}''}^{\aan})^{\rig} \rightarrow
(\LL_{\mathcal{X}''/\mathcal{X}\times_{\mathcal{Y}}\mathcal{Y}''}^{\aan})^{\rig} \rightarrow
(L\eta^{*} \LL_{\mathcal{X}\times_{\mathcal{Y}}\mathcal{Y}''/\mathcal{Y}''}^{\aan})^{\rig}[1].
\]

On the other hand, there is a natural
quasi-isomorphism $Lpr_1^* \LL_{\mathcal{X}/\mathcal{Y}}^{\aan} \cong
\LL_{\mathcal{X}\times_{\mathcal{Y}}\mathcal{Y}''/\mathcal{Y}''}^{\aan}$
by Proposition~\ref{formal-base-change}.
Hence we have $L\alpha^* \LL_{\mathcal{X}/\mathcal{Y}}^{\aan} \cong L\eta^* \LL_{\mathcal{X}\times_{\mathcal{Y}}\mathcal{Y}''/\mathcal{Y}''}^{\aan}$ by Lemma~\ref{defined} (1).
Thus taking the above triangle into account, we see that if
$(\LL_{\mathcal{X}''/\mathcal{X}\times_{\mathcal{Y}}\mathcal{Y}''}^{\aan})^{\rig}$
$=0$,
there are isomorphisms,
\[
(\LL_{\mathcal{X}/\mathcal{Y}}^{\aan})^{\rig}\cong
(L\alpha^* \LL_{\mathcal{X}/\mathcal{Y}}^{\aan})^{\rig} \cong
(L\eta^{*} \LL_{\mathcal{X}\times_{\mathcal{Y}}\mathcal{Y}''/\mathcal{Y}''}^{\aan})^{\rig} \cong
(\LL_{\mathcal{X}''/\mathcal{Y}''}^{\aan})^{\rig}.
\]
The first isomorphism in the above sequence follows from the fact that
 the functor $(\bullet)^{\rig}$
commutes with admissible blow-ups. 
Since $\eta^{\rig}:(\mathcal{X}'')^{\rig}\to (\mathcal{X}\times_{\mathcal{Y}}\mathcal{Y}'')^{\rig}$ induces the isomorphism of rigid analytic spaces,
we have 
$(\LL_{\mathcal{X}''/\mathcal{X}\times_{\mathcal{Y}}\mathcal{Y}''}^{\aan})^{\rig}=0$ from \cite[7.2.42 (1)]{GR} and \cite[Prop. 1.5]{BL3}.
\QED

We define a rigid cotangent complex $\LL_{X/Y}^{\an}$ in
$\mathsf{D}^{-}(\mathcal{O}_{X}\operatorname{-Mod})$ of $f:X \rightarrow Y$
by $(\LL^{\aan}_{\mathcal{X}/\mathcal{Y}})^{\rig}$ for some formal model
$\mathcal{X} \to \mathcal{Y}$.

For every morphism $\tilde{f}:\mathcal{X}\to \mathcal{Y}$
of admissible formal schemes, there is a natural morphism:
$\LL_{\mathcal{X}/\mathcal{Y}}^{\aan} \to \Omega_{\mathcal{X}/\mathcal{Y}}^{\aan}$
which induces an isomorphism $H_0(\LL_{\mathcal{X}/\mathcal{Y}}^{\aan})\cong \Omega_{\mathcal{X}/\mathcal{Y}}^{\aan}$ (cf. \cite[7.2.9]{GR}).
Thus we have a natural morphism
\[
\LL_{X/Y}^{\an} \to \Omega_{X/Y}^{\an}
\]
which induces an isomorphism
\[
\label{isom}
H_0(\LL_{X/Y}^{\an}) \stackrel{\sim}{\to} \Omega_{X/Y}^{\an}
\]
where $X:=\mathcal{X}^{\rig}$ and  $Y:=\mathcal{Y}^{\rig}$
and $\Omega_{X/Y}^{\an}$ is the differential module
defined in \cite[Section 1]{BL3}.

\begin{Remark}
Our construction also works in the case of
relative rigid spaces introduced and studied by
Bosch, L\"{u}tkebohmert and Raynaud (cf. \cite{BL1}
\cite{BL2} \cite{BL3} \cite{BL4}).
However, in present paper, we concentrate
our attentions on the classical case. 

In virtue of the theorem of
Raynaud~\cite[4.1]{BL1} and formal flattening theorem (See \cite[5.2]{BL2}),
as we will see below, the study
of analytic cotangent complexes of rigid analytic spaces
is reduced to the study of the formal case.
Thus Proposition~\ref{core} is important.
\end{Remark}

\begin{Proposition}
\label{rigid-base-change}
Let $f:X \to Y$ and $g:Z \to Y$ be morphisms of rigid analytic spaces.
Then there is a natural quasi-isomorphism
\[
Lpr_1^*\LL_{X/Y}^{\rig} \stackrel{\sim}{\to} \LL_{X\times_YZ/Z}^{\rig}
\]
where $pr_1$ is the first projection $X\times_YZ \to X$.
\end{Proposition}

\Proof
Take formal models $\tilde{f}:\mathcal{X}\to \mathcal{Y}$
and $\tilde{g}:\mathcal{Z} \to \mathcal{Y}$ of $f$ and $g$ respectively.
Note that the fibre product $\mathcal{X}\times_{\mathcal{Y}}\mathcal{Z}$
in the category of formal
$R$-schemes is the formal model
of $X\times_YZ$ by \cite[4.6]{BL1}.
Now our assertion follows from Proposition~\ref{formal-base-change}.
\QED

\begin{Proposition}
\label{smooth}
If $f:X \rightarrow Y$ is a smooth morphism, there is a natural
quasi-isomorphism
\[
\LL_{X/Y}^{\an} \stackrel{\sim}{\to} \Omega_{X/Y}^{\an}[0].
\]
\end{Proposition}
\Proof
Take a formal model $\tilde{f}:\XXX \to \YYY$ of
$f:X\to Y$. Then, we derive our assertion by applying
\cite[7.2.42 (1)]{GR} to $\tilde{f}:\XXX \to \YYY$.
\QED

\begin{Proposition}
\label{transitive}
Let $X \stackrel{f}{\rightarrow} Y \rightarrow Z$ be a sequence of morphisms
of rigid analytic spaces. Then, there is a natural distinguished triangle
in $\mathsf{D}^-(\mathcal{O}_{X}\operatorname{-Mod})$
\[
Lf^{*} \LL_{Y/Z}^{\an}
\rightarrow \LL_{X/Z}^{\an} \rightarrow
\LL_{X/Y}^{\an} \rightarrow
Lf^{*} \LL_{Y/Z}^{\an}[1].
\]
\end{Proposition}

\Proof
By the theorem of Raynaud \cite[4.1]{BL1}, there is  a formal model of
$X \rightarrow Y \rightarrow Z$
\[
\mathcal{X} \stackrel{\tilde{f}}{\longrightarrow} \mathcal{Y} \longrightarrow \mathcal{Z}.
\]
From the transitivity of this sequence \cite[7.2.13]{GR}, we have
\[
L\tilde{f}^{*} \LL_{\mathcal{Y}/\mathcal{Z}}^{\aan}
\rightarrow \LL_{\mathcal{X}/\mathcal{Z}}^{\aan} \rightarrow
\LL_{\mathcal{X}/\mathcal{Y}}^{\aan} \rightarrow
L\tilde{f}^{*} \LL_{\mathcal{Y}/\mathcal{Z}}^{\aan}[1].
\]
By applying the functor of derived categories
\[
\operatorname{Rig}:\mathsf{D}^{-}(\mathcal{O}_{\mathcal{X}}\operatorname{-Mod})\rightarrow
\mathsf{D}^{-}(\mathcal{O}_{X}\operatorname{-Mod}),
\]
we obtain the
triangle
\[
(L\tilde{f}^{*} \LL_{\mathcal{Y}/\mathcal{Z}}^{\aan})^{\rig}
\rightarrow (\LL_{\mathcal{X}/\mathcal{Z}}^{\aan})^{\rig} \rightarrow
(\LL_{\mathcal{X}/\mathcal{Y}}^{\aan})^{\rig} \rightarrow
(L\tilde{f}^{*} \LL_{\mathcal{Y}/\mathcal{Z}}^{\aan}[1])^{\rig}.
\]
Thus we have the required triangle by Lemma \ref{defined} (2).
\QED

\begin{Theorem}
\label{helloween}
\begin{enumerate}
\renewcommand{\labelenumi}{(\theenumi)}

  \item Let $i:Y \rightarrow X$ is a closed immersion of
rigid analytic spaces. Then the natural morphism
\[
\LL_{Y/X} \longrightarrow \LL_{Y/X}^{\an}
\]
is a quasi-isomorphism. Here $\LL_{Y/X}$ is the usual cotangent complex
associated to the morphism of ringed topoi $i:Y
\rightarrow X$.
  \item Let $f:X \rightarrow Y$ be a morphism of rigid analytic spaces.
  Then $\LL_{X/Y}^{\an}$ is a pseudo-coherent complex of $\mathcal{O}_{X}$-
  modules.
\end{enumerate}
\end{Theorem}

\Proof
(1) First of all, take a formal model
 $\tilde{i}:\mathcal{Y} \rightarrow \mathcal{X}$
of $i:Y\rightarrow X$. Note that there is a natural isomorphism
$\LL_{Y/X} \cong (\LL_{\mathcal{Y}/\mathcal{X}})^{\rig}$
by \cite[2.2.3]{IL}. Thanks to the formal flattening theorem
\cite[5.4 (b)]{BL2}, we can modify
$\tilde{i}:\mathcal{Y} \rightarrow \mathcal{X}$ to be a closed immersion.
Then the claim is reduced to \cite[7.2.10 (2)]{GR}.

(2) Let $\tilde{f}:\mathcal{X} \rightarrow \mathcal{Y}$ be a formal model
of $f:X \rightarrow Y$. Then the claim is reduced to the claim
for $\tilde{f}:\mathcal{X} \rightarrow \mathcal{Y}$ \cite[7.2.10 (1)]{GR}.
\QED

\begin{Proposition}
\label{expl}
Let $f:X \rightarrow Y$ be a closed immersion of rigid analytic spaces
and $g:Y \rightarrow Z$ a smooth morphism of rigid analytic spaces.
Let $\mathcal{I}$ be a coherent ideal of $\mathcal{O}_Y$
which defines $X$.
Then there is a quasi-isomorphism
\[
\tau_{[-1}\LL_{X/Y}^{\an}\stackrel{\sim}{\to} [0\rightarrow f^*(\mathcal{I}/\mathcal{I}^2)
\stackrel{d}{\to} f^*\Omega_{Y/Z}^{\an}\rightarrow 0].
\]
\end{Proposition}

\Proof
First, we remark that by \cite[Chapter 3, 1.2.8.1]{IL}
and Theorem~\ref{helloween} (1),
there exists a natural quasi-isomorphisms,
$\tau_{[-1}\LL_{X/Y}^{\rig}\cong \tau_{[-1}\LL_{X/Y}\cong f^*(\mathcal{I}/\mathcal{I}^2)[1]$.
Then,  we can prove our assertion in the same way as
the proof of \cite[Chapter 3, 1.2.9.1]{IL}
by using Propositions~\ref{smooth}
and \ref{transitive}.
\QED
Let $X\rightarrow Y$ be a morphism of rigid analytic spaces
and $\mathcal{S}$ a coherent $\mathcal{O}_X$-module.
Let us consider a triple $(i:X\rightarrow X', \mathcal{S},\phi)$
where $i$ is a closed immersion of rigid analytic spaces
over $Y$ with the square-zero kernel $\mathcal{I}=\Ker (\mathcal{O}_{X'}\rightarrow \mathcal{O}_{X})$
and an isomorphism of coherent $\mathcal{O}_{X}$-modules
$\phi: i^*\mathcal{I} \rightarrow \mathcal{S}$.
Then there exists a natural quasi-isomorphism
$\Ext_{\OO_X}^1(\LL_{X/X'}^{\rig},\mathcal{S})\cong \Hom_{\OO_{X}}(i^*\mathcal{I},\mathcal{S})$, since we have a natural isomorphism
$\LL_{X/X'}^{\rig}\cong i^*\mathcal{I}[1]$ by Proposition~\ref{expl}.
On the other hand, the disitinguished triangle (Proposition~\ref{transitive})
associated to the sequence
$X\rightarrow X' \rightarrow Y$ induces a homomorphism
$p:\Ext_{\mathcal{O}_{X}}^1(\LL_{X/X'}^{\an},\mathcal{S}) \rightarrow \Ext^1_{\mathcal{O}_{X}}(\LL_{X/Y}^{\an},\mathcal{S})$. Thus it gives rise to a map
\[
e:\EXan_Y(X,\mathcal{S}) \longrightarrow \Ext^1_{\mathcal{O}_{X}}(\LL_{X/Y}^{\an},\mathcal{S}),
\]
\[
(i:X\rightarrow X', \mathcal{S},\phi)\mapsto p(\phi ),
\]
where $\EXan_Y(X,\mathcal{S})$ is the set of isomorphism classes
of a triple $(i:X\rightarrow X', \mathcal{S},\phi)$.
Note that the set $\EXan_Y(X,\mathcal{S})$
is the subset of the set $\operatorname{EXal}_Y(X,\mathcal{S})$ of
isomorphism classes of the extensions of $X$ by $\mathcal{S}$
over $K$ as locally ringed spaces.
Indeed for elements $\alpha, \beta \in \EXan_Y(X,\mathcal{S})$
and an isomorphism $\phi :\alpha \stackrel{\sim}{\to}\beta $ as
locally ringed spaces over $K$ which induces the identity on $X$,
$\phi$ is actually the isomorphism as rigid analytic spaces
by Proposition~\ref{affinoid:affinoid} and \cite[6.1.3.1]{BGR}.
Hence, taking \cite[Chapter 3, 1.2.3]{IL} into account, $e$ is injective
because the composite map 
$\EXan_Y(X,\mathcal{S}) \stackrel{e}{\rightarrow} \Ext^1_{\mathcal{O}_{X}}(\LL_{X/Y}^{\an},\mathcal{S}) \to \Ext^1_{\mathcal{O}_{X}}(\LL_{X/Y},\mathcal{S})$
is equal to the composite map
$\EXan_Y(X,\mathcal{S}) \to \operatorname{EXal}_Y(X,\mathcal{S}) \stackrel{\sim}{\to} \Ext^1_{\mathcal{O}_{X}}(\LL_{X/Y},\mathcal{S})$. 

\begin{Theorem}
\label{fundamental}
The map $e$ is a bijection.
\end{Theorem}

\Proof
What we need to show is that
elements of the image $\Ext^1 (\Omega^{\rig}_{X/Y},\mathcal{S})$
in
$\Ext^1 (\Omega^1_{X/Y},\mathcal{S}) \cong \operatorname{EXal}_Y(X,\mathcal{S})$
represent the sheaves of rigid analytic spaces
which are extensions of $\mathcal{O}_{X}$ by $\mathcal{S}$.
Therefore the problem is local on $X$.
Set $X:=\SP (B)$ and $Y:=\SP (A)$.
We choose a closed immersion $i:X \cong \SP (C/I) \to Z:=\SP (C)$ 
where $C=A\lya X_1, \ldots X_r \rya $.
By Proposition~\ref{expl}, we have
\[
\Ext^1_{\mathcal{O}_X} (\LL^{\rig}_{X/Y},\mathcal{S})\cong \Hom (I/I^2,\mathcal{S})/\Image (d^*).
\]
However a map $\phi:I/I^2 \to \mathcal{S}$ 
and the canonical immersion $j:I/I^2 \to C/I^2$ define
an extension $(C/I^2\oplus \mathcal{S})/\Image ((j,\phi))^{\sim}$
of $\mathcal{O}_X$ by $\mathcal{S}$ which is a sheaf
arising from an affinoid algebra. It is easy to see that this extension
corresponds to the element $\phi$ in $\Ext^1_{\mathcal{O}_X} (\LL^{\rig}_{X/Y},\mathcal{S})$. \QED

\section{Cohomological descriptions of Local deformations}

Now by applying the results in section 2, we will prove the following theorem.

\begin{Theorem}
\label{cohomology}
Let $u:S \rightarrow S'$ be a closed immersion of
rigid analyitic spaces with the nilpotent kernel $\mathcal{I}$ of
$u^*:\OO_{S'} \rightarrow \OO_{S}$ such that
$\mathcal{I}^2=0$.
Let $f:X \rightarrow S$ be a flat morphism of rigid
analytic spaces over $K$. 

\begin{enumerate}
\renewcommand{\labelenumi}{(\theenumi)}

  \item The obstruction for the existence
  of the lifting of $f$ to $S'$ lies in
$\Ext^2_{\mathcal{O}_X}(\LL_{X/S}^{\an}, f^*\mathcal{I})$.

  \item If the obstruction $o$ is zero, the set of isomorphism classes of
the liftings of $f$ to $S'$ forms a torsor under
$\Ext^1_{\mathcal{O}_X}(\LL_{X/S}^{\an}, f^*\mathcal{I})$.

  \item Let $\tilde{f}:\tilde{X} \rightarrow S'$ be a flat 
deformation of $X$. Then,
the automorphism group of the lifting $\tilde{X}$ 
is canonically isomorphic to
$\Ext^0_{\mathcal{O}_X}(\LL_{X/S}^{\an}, f^*\mathcal{I})$.

\end{enumerate}
\end{Theorem}

\Proof
By applying Proposition~\ref{transitive}
to the sequence $X \rightarrow S \rightarrow S'$
, we have a distinguished triangle
\[
Lf^{*}\LL_{S/S'}^{\an} \rightarrow \LL_{X/S}^{\an} \rightarrow \LL_{X/S'}^{\an}
\rightarrow Lf^{*}\LL_{S/S'}^{\rig}[1].
\]
On the other hand, by Proposition~\ref{expl},
we have $\LL_{S/S'}^{\an} \cong u^* \mathcal{I}[1]$.
Thus there exists a natural isomorphism,
$\Ext^{1}_{\mathcal{O}_X}(Lf^{*}\LL_{S/S'}^{\an},\mathcal{I}\mathcal{O}_X)
\cong \Hom_{\mathcal{O}_X}(f^* \mathcal{I} \mathcal{O}_S,
\mathcal{I}\mathcal{O}_X )$.
Since $\Omega_{S/S'}^{\an}=0$,
we have $\Ext^{0}_{\mathcal{O}_X}(Lf^{*}\LL_{S/S'}^{\an},
\mathcal{I}\mathcal{O}_X)
\cong \Hom_{\mathcal{O}_X}(f^{*}\Omega_{S/S'}^{\an},
\mathcal{I}\mathcal{O}_X )$
$=0$.
Thus there is a long exact sequence
\[
\begin{CD}
0 @>>> 
 \Ext^{1}_{\mathcal{O}_X}(\LL_{X/S}^{\an},\mathcal{I}\mathcal{O}_X)
  @>>> \Ext^{1}_{\mathcal{O}_X}(\LL_{X/S'}^{\an},\mathcal{I}\mathcal{O}_X) \\
  @>>{\xi}> \Hom_{\mathcal{O}_X}(f^* \mathcal{I} \mathcal{O}_S,
  \mathcal{I}\mathcal{O}_X )  
  @>>{\delta}> \Ext^{2}_{\mathcal{O}_X}(\LL_{X/S}^{\an},\mathcal{I}\mathcal{O}_X).  \\
\end{CD}
\]
An element in  $\EXan_{S'}(X, \mathcal{I}\mathcal{O}_X)$ canonically induces
an element which lies in $\Hom_{\mathcal{O}_X}(f^* \mathcal{I} \mathcal{O}_S,
\mathcal{I}\mathcal{O}_X )$. Thus we have a map
\[
\alpha :\EXan_{S'}(X, \mathcal{I}\mathcal{O}_X) \rightarrow
\Hom_{\mathcal{O}_X}(f^* \mathcal{I} \mathcal{O}_S,
\mathcal{I}\mathcal{O}_X ).
\]
There is a natural isomorphism 
$\EXan_{S'}(X, \mathcal{I}\mathcal{O}_X) \stackrel{\sim}{\to} \Ext^{1}(\LL^{\an}_{X/S'},
\mathcal{I} \mathcal{O}_X )$ by Theorem~\ref{fundamental} and
this isomorphism identifies $\xi$ with $\alpha$.
Let $u$ be the element in $\Hom_{\mathcal{O}_X}(f^* \mathcal{I}\mathcal{O}_{S},
\mathcal{I}\mathcal{O}_{X})$
which is induced by $f$.
The existence of a flat deformation
of $X$ over $S'$ is equivalent to the existence of an element of
$\Ext^{1}_{\mathcal{O}_X}(\LL_{X/S'}^{\an},\mathcal{I}\mathcal{O}_X)$
which induces $u$ by $\xi$, so (1) follows.
Now it is clear that if $\delta u=0$, the
isomorphism classes of flat deformations is $\xi^{-1}(u)\cong
\Ext^{1}_{\mathcal{O}_{X}}(\LL_{X/S}^{\an},\mathcal{I}\mathcal{O}_{X})$.
This shows (2).
Finally the isomorphism
$H_0(\LL_{X/Y}^{\an})\cong \Omega_{X/Y}^{\an}$
implies that
$\Ext^0_{\mathcal{O}_X}(\LL_{X/S}^{\an}, \mathcal{I}\mathcal{O}_{X})
\cong \Hom(\Omega_{X/S}^{\an}, \mathcal{I}\OO_{X})$ . It is well-known that for every flat deformation of
$X$ to $S'$, its automorphism group is isomorphic to the right-hand group.
(We can show this by the completely same argument as the scheme-case (cf. \cite{SGA1}))
\QED

\section{Formal versal deformation}

\subsection{Preliminaries}
In this section, we will  consider the
deformations of rigid analytic $K$-spaces to the spectrums
of local Artin $K$-algebras.
Since every local Artin $K$-algebras is a finite $K$-vector space,
it is the topological ring whose topology are induced by
the topology of $K$. We call this topology {\em the canonical topology}.

\begin{Lemma}
\label{artin:affinoid}
Let $K$ be a complete non-Archimedean valued field.
Then, an Artin local K-algebra with the canonical topology 
is a K-affinoid algebra.
\end{Lemma}

\Proof
Since a local Artin $K$-algebra $A$ is of finite type over $K$,
we have the representation of $A$ by the quotient ring of a polynomial ring
\begin{center}
$A \cong K[X_{1},\dots ,X_{n}]/I$.
\end{center}
If we equip $A$ with canonical topology and  $K[X_{1},\dots ,X_{n}]/I$
with Gauss norm, this is an isomorphism of complete topological rings.
 The Tate algebra $K \lya X_1, \dots ,X_n \rya $ is flat over 
$K[X_1, \dots , X_n ]$. 
Thus there exists a canonical isomorphism
\begin{center}
$K[X_{1},\dots ,X_{n}]/I \cong K \lya X_1, \dots ,X_n \rya /IK \lya X_1, \dots ,X_n \rya$.
\end{center}
This completes the proof.
\QED

\begin{Remark}
By the above lemma,
 we can attach an affinoid space $\SP (A)$ to any Artin
local $K$-algebra $A$.
Since every homomorphism of affinoid algebras as $K$-algebras
is automatically a continuous homomorphism,
we can view the category of Artin local $K$-algebras
as the full subcategory of the category of $K$-affinoid algebras.

\end{Remark}
\begin{Definition}
For a rigid analytic space $X$ over $K$,
a {\em local deformation functor} $\mathbb{D}_X$ of $X$ 
 is the functor defined as follows:
\begin{center}
$\mathbb{D}_{X}$ $:$ $\biggl($ Artin local $K$-algebras 
with residue field $K$ $\biggr)$
$\longrightarrow$ $\biggl($ Sets $\biggr)$ 
\end{center}
Let $A$ be an Artin local $K$-algebra with the maximal ideal
 $\mathfrak{m}_A$ and the residue field $K$.
Then the set $\mathbb{D}_{X}(A)$ is the isomorphism classes
of pairs 
$(f:\tilde{X}\rightarrow \SP (A),\alpha :\tilde{X}\times_{A}A/\mathfrak{m}_A$
$ \xrightarrow{\thicksim} X )$
of a flat morphism $f$
of rigid analytic spaces over $K$ and
$\alpha$ is an isomorphism of rigid analytic spaces over $K$.
\end{Definition}

\subsection{Schlessinger's theory}

\medskip

Let $X$ be a rigid analytic spaces over $\SP (K)$.

\begin{Definition}
Let $\OO$ be a complete local noetherian $K$-algebra with
maximal ideal $\mathfrak{m}$ with residue field $K$ .
Let $\{ X_n \}_{n \ge 0}$ be a family 
of deformations of $X$ to $\SP (\OO/\mathfrak{m}^{n+1})$
such that $X_n$ is a flat deformation to $\SP (\OO/\mathfrak{m}^{n+1})$
and $X_n\times_{\SP (\OO/\mathfrak{m}^{n+1})}\SP (\OO/\mathfrak{m}^{m+1})$
$\cong X_m$ for $m \le n$. We say that a pair $(\OO ,\{X_n\}_{n \ge 0})$ is 
{\em a formal versal deformation} of $X$ if it satisfies
the following conditions.

\medskip

\begin{enumerate}
\renewcommand{\labelenumi}{(\theenumi)}

\item Suppose that $A$ is an Artin local $K$-algebra with the
maximal ideal $\mathfrak{m}_A$
such that $\mathfrak{m}_A^{n+1}=0$.
If $\XX$ is a deformation of $X$ to $\SP (A)$,
there exists a local homomorphism
$f:\OO/\mathfrak{m}^{n+1} \rightarrow A$ such that $\XX$ is isomorphic to
$X_n\times_{\SP (\OO/\mathfrak{m}^{n+1})}\SP (A).$

\medskip

\item If $A=K[\epsilon]/(\epsilon^2)$, such $f$ is unique.

\end{enumerate}
\end{Definition}

The existence of a formal versal deformation of $X$
is equivalent to the existence of a prorepresental hull
of the functor $\mathbb{D}_{X}$,
in the sense of Schlessinger (cf. \cite{Sch}).

\medskip

We have a convenient criterion for the existence of a prorepresental
hull.

\begin{Theorem}[Schlessinger \cite{Sch}]
Let $A' \rightarrow A$ and $A'' \rightarrow A$ be morphisms
of Artin local $K$-algebras. Consider the natural map
\begin{center}
$F:\mathbb{D}_X(A'\times_AA'') \rightarrow $
$\mathbb{D}_X(A')\times_{\mathbb{D}_X(A)}\mathbb{D}_X(A'')$.
\end{center}
The $\mathbb{D}_{X}$ has a prorepresental hull if and only if
the following conditions are satisfied:

\textup{(H1)} If $A'' \rightarrow A$ is a small surjection, $F$ is surjective.

\textup{(H2)} $F$ is bijective when $A=K$, $A''=K[\epsilon]/(\epsilon^2)$.

\textup{(H3)} $\DIM_K\mathbb{D}_X(K[\epsilon]/(\epsilon^2)<\infty$.
\end{Theorem}

\subsection{Existence of a formal versal deformation}

\begin{Theorem}
\label{formal:versal}
Let $X$ be a proper rigid analytic space over $K$.
Then, there exists a formal versal deformation of $X$.
\end{Theorem}

First we will show (H1).

\begin{Claim}

The functor $\mathbb{D}_X$ satisfies \textup{(H1)}.
\end{Claim}

\Proof
Let $(\xi',\xi'')$ be an element of
$\mathbb{D}_X(A')\times_{\mathbb{D}_X(A)}\mathbb{D}_{X}(A'')$.
We put
\begin{align*}
\xi&:=(\XX /\SP (A),\phi :\XX\times_AK \cong X),
\\
\xi' &:=(\XX' /\SP (A'),\phi': \XX'\times_AK \cong X), \\
\xi'' &:=(\XX'' /\SP (A''),\phi: \XX''\times_AK \cong X),
\end{align*}
such that $\alpha^*\xi'=\xi$ and $\beta^*\xi''=\xi$ where
$A' \xrightarrow{\alpha} A$, $A'' \xrightarrow{\beta} A$.
Fix an affinoid open set $\SP (R)$ of $\XX$.
Note that $|X|=|\XX|=|\XX'|=|\XX''|$ . 
Here $|\bullet|$ means the underlying set.
The subspace $\XX'|_{\SP (R)}$ (resp. $\XX''|_{\SP (R)}$) of $\XX'$ (resp. $\XX''$) which is
the lifting of $\SP (R)$ is an affinoid space by Proposition~\ref{affinoid:affinoid}.
We set $\SP (R')=\XX'|_{\SP (R)}$,
$\SP (R'')=\XX''|_{\SP (R)}$ and
$s:R'\rightarrow R$, $t:R''\rightarrow R$.
To show our claim, it suffices to prove that $R'\times_RR''$
is an affinoid algebra.
First we define a norm on $R'\times R''$ by the direct sum of
the norms of $R'$ and $R''$. This norm induces the
topology on $R'\times_RR''$. This topology is complete.
Indeed, if $\{(a_n,b_n)\}_{n\le 1}$ is a Cauchy sequence of $R'\times_RR''$,
$\{s(a_n)=t(b_n)\}_{n\le 1}$ is also a Cauchy sequence of $R$ since
$|t(b_n)-t(b_m)| \le |b_n-b_m|$.
Now what to do is to construct a continuous 
surjective map form an affinoid algebra to $R'\times_RR''$.
To this end, take a generator $\{\xi_1,\ldots,\xi_r\}$ of $\Ker (t)$
such that $|\xi_i|_{sp} \le 1$ for $1\le i \le r$.
Next put $R'=K \lya X_1,\ldots,X_n \rya /I$. We can choose elements
$\eta_1,\ldots,\eta_n$ in $R''$ such that
$s(X_i)=t(\eta_i)$ and
$|\eta_i|_{sp} \le 1$ for $1\le i \le n$.
Indeed, since $t$ is a small surjection,
$|r|_{sp}=|t(r)|_{sp}$ for $r\in R''$
and $|s(X_i)|_{sp} \le |X_i|_{sp} \le |X_i|\le 1$ for $1\le i \le n$.
Therefore, we have the following homomorphism
\begin{center}
$K[S_1, \dots,S_n,T_1,\dots,T_r] \rightarrow R'\times_RR''$
\end{center}
defined by $S_i\rightarrow (X_i,\eta_i)$ and $T_j \rightarrow (0,\xi_j)$.
Since $|X_i|_{sp}\le 1$, $|\eta_i|_{sp}\le 1$ and $|\xi_j|_{sp}\le 1$
for all $i$ and $j$, this homomorphism is uniquely extended to
the continuous homomorphism (See \cite[3.4.7]{FV})
\begin{center}
$K \lya S_1,\ldots,S_n,T_1,\cdots,T_r \rya \rightarrow R'\times_RR''$.
\end{center}
Furthermore, from the construction, it is clear that
this continuous homomorphism 
is surjective. Hence we see (H1).
\QED

Next we show (H2) and (H3).

\begin{Claim}
The functor $\mathbb{D}_X$ satisfies \textup{(H2)}.
\end{Claim}

\Proof
We can show this by the completely same argument as the scheme-case (cf. \cite{Sch}).
\QED

\begin{Claim}
The functor $\mathbb{D}_X$ satisfies \textup{(H3)}.
\end{Claim}

\Proof
Let $f:X \rightarrow \SP (K)$ be a structure morphism.
Then, the isomorphism of derived functors
\[
Rf_{*}R\mathcal{H}om(-,\mathcal{O}_{X}) \cong R\Hom(-,\mathcal{O}_{X})
\]
induces the spectral sequence
\[
E^{p,q}_{2}=R^{q}f_{*}\mathcal{E}xt^{p}(\LL_{X/K}^{\an},\mathcal{O}_{X})
\Rightarrow 
\Ext^{p+q}(\LL^{\an}_{X/K},\mathcal{O}_{X}).
\]
Then it suffices to show that $R^qf_{*}\mathcal{E}xt^p(\LL_{X/K}^{\an},\mathcal{O}_{X})$
is coherent for all $p$ and $q$.
Note that by Theorem~\ref{helloween} (2), $\LL_{X/K}^{\an}$ is pseudo-coherent.
Thus $R\mathcal{H}om(\LL_{X/K}^{\an},\mathcal{O}_{X})$ is pseudo-coherent by 
\cite[Chapter 0, 12.3.3]{EGA3} and \cite[Chapter 1, 7.3]{Har}.
By Kiehl's finiteness theorem,
$R^qf_{*}\mathcal{E}xt^p(\LL_{X/K}^{\an},\mathcal{O}_{X})$ is a finite $K$-vector space for all $p$ and $q$.
\QED

By the above three claims,
 we can complete the proof of Theorem~\ref{formal:versal} by Schlessinger's
 criterion.

\section{Towards the global moduli theory via Rigid geometry}

In Section 4, we proved the existence
of a versal family in the formal sense for deformations of rigid analytic spaces. 
However the analogy of rigid geometry with complex analytic geometry
gives us a deeper question.
When one compares the deformation theory of complex analytic spaces
with one of the algebraic categories,
no one doubts that the most important
advantage is the existence of
versal family of deformations proven by
Kuranishi and Grauert (cf. \cite{Gra} \cite{Kura}).
Formal deformations of algebraic varieties are
not necessarily algebraic.
Thus there is no algebraic analogue of
the theorem of Kuranishi and Grauert. As we know, there is no logical relation
between complex geometry and rigid geometry.
However one can conjecture the following.

\begin{Conjecture}
Let $X$ be a proper rigid analytic spaces over $K$.
Then there exists
a flat morphism of rigid analytic spaces
$F: \mathfrak{X} \longrightarrow \mathfrak{S}$
and a $K$-rational point $p$ of $\mathfrak{S}$
such that the completion of $F$ at $p$ is isomorphic to
the formal versal deformation of $X$.
\end{Conjecture}

This assertion can be viewed as a fairly precise non-Archimedean analogue
of the existence theorem of versal families for deformations of
complex analytic spaces due to Kuranishi and Grauert.
Let $\mathcal{S}$ be a rigid analytic torus, i.e.,
a quotient space $T^{\aan}/\Gamma$ where
$T$ is a split $K$-torus and $\Gamma$
is a torsion-free lattice of rank $\operatorname{dim}T$.
Then, we can prove that the conjecture holds for $\mathcal{S}$
by using p-adic uniformization theory due to
Bosch-L\"{u}tkebohmert-Raynaud.
Unfortunately, at the time of writting this paper, the author
do not have a proof of this conjecture for general rigid analytic spaces.
But the author expects that this conjecture is true and propose it.

Let us give one sufficient condition which implies the conjecture.
Let $K$ be a discrete valuation field and $R$ its ring of integers
with residue field $k$. Let $X$ be a proper rigid analytic
$K$-space.
To prove the conjecture for $X$, it suffices to show the existence 
of a locally noetherian adic formal scheme over $\Spf R$ which
satisfies the followings (cf. \cite{Be} section 0.2).
\begin{enumerate}

 \item Its reduction is a scheme locally of finite type over
$\Spec k$.
\item It is a formal model of the formal deformation family of a rigid
analytic spaces $X$.

\end{enumerate}

However it seems difficult to prove the conjecture in general.

Let us explain why this problem is important.
Suppose that
we want to construct a moduli space of interesting geometric objects.
From the stack theoretic viewpoint,
versal spaces for deformations of them
(here ``versal space" is not in formal sense
but has a geometric structure such as a scheme
(resp. complex analytic or rigid analytic etc...))
are local components of the smooth cover of an algebraic
(resp. complex analytic, rigid analytic) stack (cf. \cite{ArtVer}).
Thus, roughly speaking, this conjecture says
that in rigid geometry, the existence of local deformation theory
implies the global moduli stack which is represented by a rigid analytic stack.

Furthermore this conjecture will
be important to the
geometry of algebraic schemes
over an arbitrary field.
Let $X$ be an algebraic scheme over an arbitrary  field $k$.
By considering the constant deformation 
$\hat{X} \rightarrow \Spf (k\lformal T \rformal )$.
one can associate a rigid analytic space
$\hat{X}^{rig} \rightarrow \Spf (k\lformal T \rformal)^{rig}.$
By this technique the deformation of $X$ to a family of algebraic schemes
is contained in the rigid analytic versal family
if it exists. Therefore the above conjecture 
is of prime importance not only in rigid geometry
but also in algebraic geometry.

\renewcommand{\thesection}{}

\renewcommand{\theTheorem}{A.\arabic{Theorem}}

\renewcommand{\theClaim}{A.\arabic{Theorem}.\arabic{Claim}}

\renewcommand{\theequation}{A.\arabic{Theorem}.\arabic{Claim}}
\setcounter{Theorem}{0}

\section*{Appendix}

In this appendix, we prove a convenient criterion of the
existence of a non-trivial square-zero extension of a ringed topos.

Let $A$ be a ring and $B$ an $A$-algebra.
Let $M$ be a $B$-module and $\LL_{B/A}$
the cotangent complex of the structure homomorphism $A \to B$.
By $\operatorname{EXal}_A(B,M)$, we denote the set of
isomorphism classes of sequare-zero extension of $B$
by $M$ over $A$.
By the fundamental theorem due to L. Illusie (\cite[Chaptor 3, 1.2.3]{IL}),
there exists a natural bijection
\[
\phi:\Ext_A^1(\LL_{B/A},M) \stackrel{\sim}{\to} \operatorname{EXal}_A(B,M).
\]
On the other hand, there exists a natural
homomorphism 
\[
\pi:\Ext_A^1(\Omega^1_{B/A},M) \to \Ext_A^1(\LL_{B/A},M)
\]
which is induced by $\LL_{B/A} \to \Omega^1_{B/A}$ (cf. \cite[Chapter 2, 1.2.4]{IL}).
\begin{Theorem}
\label{inj}
The natural map $\pi:\Ext_A^1(\Omega^1_{B/A},M) \to \Ext_A^1(\LL_{B/A},M)$
is injective. 
In particular, if $\Ext_A^1(\Omega^1_{B/A},M) \ne 0$, there exists a
non-trivial
square-zero extension of $B$ by $M$ over $A$.

\end{Theorem}

\Proof
First, note that
it suffices to show that the map $\psi:=\phi \circ \pi:\Ext_A^1(\Omega^1_{B/A},M) \to \operatorname{EXal}_A(B,M)$ is injective.
Let us construct explicitly the map $\psi$ (cf. \cite[Chapter 3, 1.1.8]{IL}). For an element $\xi$ in $\Ext_A^1(\Omega^1_{B/A},M)$,
let $(0 \to M \stackrel{\alpha}{\to} N \stackrel{\beta}{\to} \Omega^1_{B/A} \to 0)$ be the corresponding short
exact sequence. We define an $A$-algebra structure of $B\oplus N$
by $(b,n)+(b',n'):=(b+b',n+n')$ and $(b,n)\cdot(b',n'):=(b\cdot b',b\cdot n'+b'\cdot n)$. We also define an $A$-algebra structure of $B\oplus \Omega^1_{B/A}$
by the same way. Consider the following diagram
\[
 \xymatrix{
 0 \ar[r]  & M \ar[r] \ar[d]^{\operatorname{Id}} & C \ar[r]^{\operatorname{pr}_2} \ar[d]^{\operatorname{pr}_1} & B \ar[r] \ar[d]^{(\operatorname{Id},d_{B/A})} & 0 \\
 0 \ar[r] & M \ar[r]_{(0,\alpha)} & B\oplus N \ar[r]_{(\operatorname{Id},\beta)} & B\oplus \Omega^1_{B/A} \ar[r] & 0, \\
  }
\]
where $C$ is the fibre product $(B\oplus N) \times_{(B\oplus \Omega^1_{B/A})}B$.
Then we define $\psi (\xi)$ by $(0 \to M \to C \to B \to 0)$.
Since $\operatorname{EXal}_A(B,M)$ has a group structure by $\phi$,
it suffices to prove the following claim.

\begin{Claim}
Under the same assumption as above, if the exact sequence
$(0 \to M \to C \to B \to 0)$ has a splitting $a:B \to C$,
the exact sequence $(0 \to M \to N \to \Omega^1_{B/A} \to 0)$
also has a splitting $\Omega^1_{B/A} \to N$.
\end{Claim}

\Proof
Let $\operatorname{Der}_A(B,N)$ be the set of $A$-derivations of $B$ to $N$.
Then note that we have a natural isomorphisms 
$\Hom_B(\Omega^1_{B/A},N) \stackrel{\sim}{\to} \operatorname{Der}_A(B,N)$
and
$\operatorname{Der}_A(B,N)\stackrel{\sim}{\to} \Hom_{A\operatorname{-alg/B}}(B,B\oplus N)$
(cf. \cite[Chapter 2, 1.1.1.4 and 1.1.2.6]{IL}).
Let $\delta:\Omega^1_{B/A} \to N$ be a homomorphism  that corresponds to the element $\operatorname{pr}_1\circ a$ in $\Hom_{A\operatorname{-alg/B}}(B,B\oplus N)$ by the above isomorphisms.
Then it is easy to see that $\delta$ is a splitting of
$\beta: N \to \Omega^1_{B/A} \to 0$.\QED
Thus we completes the proof of the theorem. \QED
Let $f:X \to Y$ be a morphism of ringed topoi.
Let $\Omega^1_{X/Y}$ be a K\"{a}hlar differential module of $f^{-1}(\mathcal{O}_{Y}) \to \mathcal{O}_X$, and $\LL_{X/Y}$ the cotangent complex (cf. \cite[chapter 2]{IL}).
By  $\operatorname{EXal}_Y(X,\mathcal{M})$, we denote the set of
isomorphism classes of sequare-zero extension of $X$
by a $\mathcal{O}_X$-module $\mathcal{M}$ over $Y$.
By the same procedure with the above,
we have a natural bijection
$\phi:\Ext_Y^1(\LL_{X/Y},\mathcal{M}) \stackrel{\sim}{\to} \operatorname{EXal}_Y(X,\mathcal{M})$
and a homomorphism
$\pi:\Ext_Y^1(\Omega^1_{X/Y},\mathcal{M}) \to \Ext_Y^1(\LL_{X/Y},\mathcal{M})$.
\begin{Corollary}
Let $X \to Y$ be a morphism of ringed topoi and $\mathcal{M}$ a $\mathcal{O}_X$-module.
The natural map $\pi:\Ext_Y^1(\Omega^1_{X/Y},\mathcal{M}) \to \Ext_A^1(\LL_{X/Y},\mathcal{M})$
is injective. 
In particular, if $\Ext_Y^1(\Omega^1_{X/Y},\mathcal{M}) \ne 0$, there exists a
non-trivial
square-zero extension of $X$ by $\mathcal{M}$ over $Y$.

\end{Corollary}


\begin{thebibliography}{99}



\bibitem{ArtVer}
M. Artin,
Versal deformations and Algebraic stacks,
Invent. Math. vol. 27 (1974) 165--184

\bibitem{Be}
P. Berthelot,
Cohomologie rigide et cohomologie rigide support propres,
Pr\'{e}publication IRMAR 96-03, 89 pages (1996).





\bibitem{BGR}
S. Bosch, U. Guntzer and R. Remmert,
Non-Archimedean analysis.
A systematic approach to rigid analytic geometry.
Grundlehren der math. Wiss., 261 Springer-Verlag, Berlin, 1984.

\bibitem{BL}
S. Bosch and W. L\"{u}tkebohmert,
Degenerating abelian varieties,
Topology, vol. 30 No. 4, (1991) 653--698

\bibitem{BL1}
S. Bosch and W. L\"{u}tkebohmert.
Formal and Rigid Geometry I.
Math. Ann, 295 (1993), no.2, 291--317.

\bibitem{BL2}
S. Bosch and W. L\"{u}tkebohmert.
Formal and Rigid Geometry II.
Math. Ann, 296 (1993), no.3, 403--429.

\bibitem{BL3}
S. Bosch, W. L\"{u}tkebohmert and M. Raynaud.
Formal and Rigid Geometry III.
Math. Ann, 302 (1995), no.1, 1--29.

\bibitem{BL4}
S. Bosch, W. L\"{u}tkebohmert and M. Raynaud.
Formal and Rigid Geometry IV.
Invent. Math., 119 (1995), no.2, 361--398.

\bibitem{Bou}
N. Bourbaki,
El\'{e}ment de Math\'{e}matique:
Alg\'{e}bre Commutative,
Pris: Hermann 1961--1965


\bibitem{EGA1}
J. Dioudonn\'{e} and A. Grothendieck,
\'Elements de G\'{e}om\'{e}trie Alg\'{e}briques I,
Publ. Math. IHES 4 (1960)


\bibitem{EGA3}
J. Dioudonn\'{e} and A. Grothendieck,
\'Elements de G\'{e}om\'{e}trie Alg\'{e}briques Chapitre III,
Publ. Math. IHES 11 (1961)

\bibitem{FV}
J. Fresnel and M. van der Put.
Rigid Analytic Geometry and Its Applications.
Progress in Math. vol.218 (2003)

\bibitem{GR}
O. Gabber and L. Ramero,
Almost ring theory,
Lecture Note in Math. 1800 (2003)

\bibitem{SGA1}
A. Grothendieck,
Rev\^etment \'Etale et Groupe Fondamental,
Lecture Note in Math. 224 Springer-Verlag, Berlin, 1971


\bibitem{Gra}
H. Grauert, Der Satz von Kuranishi f\"{u}r komplexe R\"{a}ume,
Invent. Math. 25 (1974) 107--142



\bibitem{Har}
R. Harshorne,
Residues and duality.
Lecture Note in Math. 20 (1966)

\bibitem{IL}
L. Illusie,
Complexe Cotangent et D\'{e}formations 1
Lecture Note in Math. 239 (1971)

\bibitem{KS}
M. Kashiwara and P. Scapira,
Sheaves on manifolds,
Die Grundlehren der Math. Wiss., 292 (1990) Springer




\bibitem{Kura}
M. Kuranishi,
On the locally complete families of complex analytic structures,
Ann. Math. (2) vol. 75, No. 3 (1962) 536--577

\bibitem{RZ}
M. Rapoport and T. Zink,
Period spaces for p-divisible groups.
Annals of Math. Studies. Princeton University Press. (1996)

\bibitem{Sch}
M. Schlessinger,
Functor of Artin rings,
Trans. Amer. Soc. 130 (1968) 208--222






\end{thebibliography}
\end{document}